\documentclass[leqno]{article}
\usepackage[latin1]{inputenc}
\usepackage{hyperref}
\usepackage{enumerate}
\usepackage{amsmath}
\usepackage{amssymb}
\usepackage{amsthm}                                         
\usepackage{url}
\usepackage{esvect}
\usepackage{graphicx}
\usepackage{epsfig}
\usepackage[dvipsnames]{xcolor}
\newtheorem{theorem}{\bf Theorem}

\newtheorem{lemma}[theorem]{\bf Lemma}

\newtheorem{proposition}[theorem]{\bf Proposition}

\newtheorem{remark}[theorem]{\bf Remark}

\numberwithin{equation}{section}
\numberwithin{theorem}{section}
\numberwithin{figure}{section}

\def\R{\mathbb{R}}
\def\L{\mathbb{L}}
\def\R{\mathbb{R}}

\newcommand{\p}{\shortmid}
\newcommand{\pp}{\shortparallel}
\def\tu{\tilde{u}}
\def\tx{\tilde{x}}
\def\ty{\tilde{y}}
\def\wp{{\widetilde \psi}}
\def\wO{{\widetilde \Omega}}
\def\wph{{\widetilde \phi}}
\def\wW{{\widetilde W}}

\makeatletter
\newsavebox{\@brx}
\newcommand{\llangle}[1][]{\savebox{\@brx}{\(\m@th{#1\langle}\)}%
  \mathopen{\copy\@brx\kern-0.5\wd\@brx\usebox{\@brx}}}
\newcommand{\rrangle}[1][]{\savebox{\@brx}{\(\m@th{#1\rangle}\)}%
  \mathclose{\copy\@brx\kern-0.5\wd\@brx\usebox{\@brx}}}
\makeatother

\begin{document}
\renewcommand{\thefootnote}{}
\footnotetext{Research partially supported by Ministerio de Econom\'ia y Competitividad Grant No: MTM2016-80313-P,  Junta de Andaluc\'ia Grant No. A-FQM-139-UGR18 and the `Maria de Maeztu'' Excellence
Unit IMAG, reference CEX2020-001105-M, funded by
MCIN/AEI/10.13039/501100011033} 

\title{Geometry of  $[\varphi,\vec{e}_{3}]$-minimal surfaces in $\mathbb{R}^{3}$}
\author{Antonio 
Mart\'inez and A. L. Mart\'inez-Trivi\~no}
\date{}
\vspace{.1in}
\maketitle
{\small
\noindent Department of Geometry and Topology, University of Granada, E-18071 Granada, Spain\\ \\
e-mails: amartine@ugr.es, aluismartinez@ugr.es 


}
\everymath={\displaystyle}
\begin{abstract}
 In this survey we report a general and systematic approach to study $[\varphi,\vec{e}_{3}]$-minimal surfaces in $\mathbb{R}^{3}$ from a geometric viewpoint and show some  fundamental results obtained in the recent development of this theory.
 \end{abstract}
 \tableofcontents
\section{Introduction}
In this survey, we would like to report some recent results about critical points of the weighted area functional
\begin{equation}
\label{areafunctional}
\mathcal{A}^{\varphi}(\Sigma)=\int_{\Sigma}\, e^{\varphi}\, d\Sigma,
\end{equation}
on surfaces $\Sigma$ in a domain $\mathfrak{D}^3\subset \mathbb{R}^3$ when $\varphi$ is the restriction on $\Sigma$ of a smooth function depending only on the last coordinate of  $\mathfrak{D}^3$ and where $d\Sigma$ denotes the volume element induced by the Euclidean metric $\langle\cdot,\cdot\rangle$  in $\mathbb{R}^3$. 

\

The Euler-Lagrange equation of \eqref{areafunctional} is given in terms of the mean curvature vector $\textbf{H}$ of $\Sigma$ as follows
\begin{equation}
\label{meancurvature}
\textbf{H}=(\overline{\nabla}\varphi)^{\perp} = \dot{\varphi} \ \vec{e}_3^{\,\perp},
\end{equation}
where $\perp$ denotes the projection on the normal bundle, $\overline{\nabla}$ stands the usual gradient operator in $\mathbb{R}^{3}$ and  $(\ \dot{ }\ )$  denotes derivate with respect to the third coordinate.   Ilmanen in \cite{Ilm94},  proved that \eqref{meancurvature} means also that $\Sigma$ is a \textit{minimal} immersion in the so called Ilmanen's space, $(\mathfrak{D}^3, g_\varphi)$ that is, $\mathfrak{D}^3$ endowed  with the conformally  Euclidean changed  metric \begin{equation}
\label{metricIlmanen}
g_\varphi:=e^{\varphi}\langle\cdot,\cdot\rangle.
\end{equation}

\

From a Physical point of view, see \cite[pp. 173-187]{Poisson}, the equation \eqref{meancurvature} gives the equilibrium condition of a flexible inextensible surface in the absence of intrinsic forces under the gravitational force field $${\cal F}:= \overline{\nabla}\mathrm{e}^\varphi=(0,0,\dot{\varphi}\,\mathrm{e}^\varphi).$$

\

Any surface satisfying $\eqref{meancurvature}$ will be called  {\sl $[\varphi, \vec{e}_3]$-minimal } 
and if  $\Sigma$ is   the vertical graph of a function $u:\Omega\subseteq \R^2 \longrightarrow \R$, we  also refer to $u$ as  $[\varphi, \vec{e}_3]$-minimal. Hence, $u$ is  $[\varphi, \vec{e}_3]$-minimal if and only if it solves the so called $[\varphi, \vec{e}_3]$-minimal equation, \begin{equation}
(1+u_{x}^2)u_{yy} + (1+u_{y}^2)u_{xx} - 2u_{y}u_{x}u_{xy}= \dot{\varphi}(u)W^2, \quad (x,y)\in\Omega, \label{fe}
\end{equation} 
where $W=\sqrt{1+u_{x}^{2}+u_{y}^{2}}$. If $\Omega$ is simply connected then, the  Poincare's Lemma gives that the equation \eqref{fe} is equivalent to integrability of the following differential system,
\begin{equation}
\label{intsystem}
\phi_{xx}=e^{\varphi(u)}\frac{1+u_{x}^{2}}{W}, \ \ \phi_{xy}=e^{\varphi(u)}\frac{u_{x}u_{y}}{W}, \ \ \phi_{yy}=e^{\varphi(u)}\frac{1+u_{y}^{2}}{W},
\end{equation}
 for a convex function $\phi:\Omega\rightarrow\mathbb{R}$ unique, up to linear polynomials.
 
 \
 
This kind of surfaces has been widely studied specially from the viewpoint of calculus of variations. Classical results about the Euler equation and the existence and regularity for  the solutions of the Plateau problem for \eqref{areafunctional} can be found in \cite{BHT, H1, H2, HK, T}.
But contributions from a more geometric viewpoint only has been given for some  particular  functions $\varphi$, namely,\begin{itemize}
\item for translating solitons:  if $\varphi$ is just the height function, $\varphi(z)=z$, that is,  surfaces such that $$ t \rightarrow \Sigma + t \vec{e}_3 $$
is a mean curvature flow, i.e. the normal component of the  velocity at  each point is equal to the mean curvature at that point. Recent advances in the understanding of their local and global geometry can be found in \cite{CSS, HMW, HIMW, HIMW2,MSHS1, MSHS2, SX, W}
\item for singular $\alpha$-minimal surfaces: if $\varphi(z) =\alpha \log z$, $z>0$, $\alpha$=const. (when $\alpha=-2$,  $\Sigma$ is a minimal surface in the Poincaré upper half hyperbolic space model and when  $\alpha=1$,  $\Sigma$ describes the shape of a ``hanging roof'', i.e. a heavy surface in a gravitational field that, according to the architect F. Otto \cite[p. 290]{Otto} are of importance for the construction of perfect domes). We refer to \cite{BHT,D, DH, Du, Rafa1, Rafa2, Rafa3, Rafa4, MM, N} for some  progress  in this family.
\end{itemize}
The aim of this paper is to develop a systematic and general approach to the study of  $[\varphi, \vec{e}_3]$-minimal surfaces from a geometric point of view. Since this class of surfaces is too large, actually, we could find almost any geometric asymptotic  behavior, it will be necessary to give   some conditions on the function $\varphi$. Here, as a general assumption we will always consider $\varphi$ strictly monotone, that is,
\begin{align}\label{hip}
& \varphi:]a,b[\subseteq \R \rightarrow \mathbb{R} \text{ is a strictly increasing (or decreasing) function}\\
& \text{and  $\Sigma\subset \mathfrak{D}^3=\R^2 \times ]a,b[$.}\nonumber
\end{align}

\section{The most symmetric examples}
As $\varphi$  is considered so arbitrary, we are only going to describe $[\varphi,\vec{e}_{3}]$-minimal surfaces invariant by either horizontal translations or vertical rotations.
\subsubsection*{{\sc $\bullet$ The one-dimensional variational problem}}
Let us consider the one-dimensional variational integral
\begin{equation} \label{oned}
{\cal I}_\varphi(u) = \int_I \mathrm{e}^\varphi \sqrt{1 + u'(x)^2} dx,
\end{equation}
with $u:I\subseteq \mathbb{R}\rightarrow ]a,b[$ a differentiable function. It is easy to check that extremals for ${\cal I}_\varphi$ must be solutions of the following ODE,
%
%
\begin{equation}
\label{grafollano}
u''(x)=\dot{\varphi}(u)(1+u'(x)^{2}).
\end{equation}
Looking for complete examples, we will assume that $\varphi:]a,+\infty[\rightarrow\mathbb{R}$, $a\in\mathbb{R}\cup\{-\infty\}$ is strictly monotone. Then,  taking  $z=\varphi(u)$ and $u'=\text{tan}(v)$, we get that \eqref{grafollano} is equivalent to
\begin{equation}
\label{sistema}
\left.
\begin{array}{rcl}
     v' & = & \lambda(z),
  \\ z' & = & \lambda(z)\tan(v)

\end{array}
\right\}
\end{equation}
where $\lambda$ is the function defined by $\lambda(z)=\dot{\varphi}(\varphi^{-1}(z))$ on $\varphi(]a,+\infty[)$.

\begin{figure}[htb]
\label{dfases}
\begin{center}
\includegraphics[width=.4\textwidth]{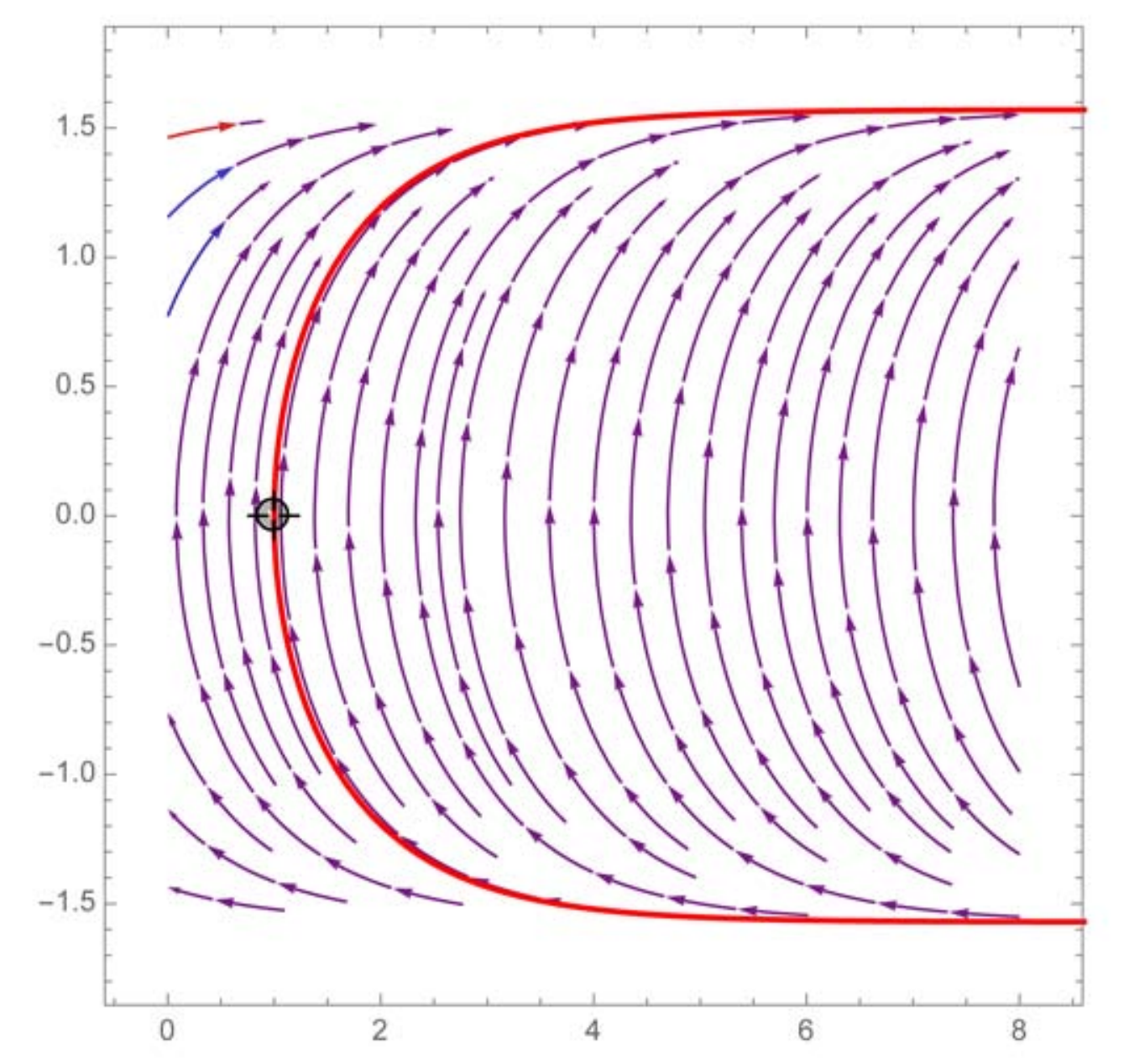}
\caption{Phase portrait of \eqref{sistema}.}
\end{center}
\end{figure}
Thus,  $e^{z}\text{cos}(v)$ is constant along the solutions of \eqref{sistema} and from the Phase portrait of \eqref{sistema} (see Figure \ref{dfases}), for any maximal solution $u$ of \eqref{grafollano}, there exists a unique $x_{0}\in I$ such that $v(x_{0})=0$. It is not a restriction to consider  $x_{0}=0$. and   $u$ satisfying
\begin{equation}
\label{initialcond}
u(0)=u_0, \ \ \ \ u'(0)=0.
\end{equation}
In this case,  the solution to \eqref{grafollano}-\eqref{initialcond},  is given by
\begin{equation}
 \label{solution}
u(x)=(\mathcal{X}\circ\varphi)^{-1}(x), \text{ with } \ \ \mathcal{X}(z)=\int_{z_{0}}^{z}\frac{d\tau}{\vert h(\tau)\vert\sqrt{e^{2(\tau-z_{0})}-1}},
\end{equation}
where $z_{0}=\varphi(u_{0})$. Thus,  we obtain that $u$ is even and it is defined in the interval $]-\Lambda_{u_{0}},\Lambda_{u_{0}}[$, where
\begin{equation}
\label{Lambda}
\Lambda_{u_{0}}=\lim_{u\rightarrow+\infty}\int_{\varphi(u_{0})}^{\varphi(u)}\frac{d\tau}{\vert h(\tau)\vert\sqrt{e^{2(\tau-z_{0})}-1}}.
\end{equation}
%
\begin{theorem}[{\rm \cite[Theorem 3.2 and Theorem 3.3]{MM}}]\label{t2} Let $\varphi:\ ]a,+\infty[\ \rightarrow \ ]b,c[$,  $ a,b\in \R\cup\{-\infty\}$,  $c\in  \R\cup\{+\infty\}$ be a strictly increasing diffeomorphism. Then the solution $u$  of \eqref{grafollano}-\eqref{initialcond} is defined in $]-\Lambda_{u_{0}},\Lambda_{u_{0}}[$, $\Lambda_{u_{0}} \in \R^{+}\cup\{+\infty\}$, it is convex,  symmetric about the $y$-axis and has a minimum at $x=0$. Moreover, 
\begin{enumerate}
\item if $c<+\infty$, then $\Lambda_{u_{0}}=+\infty$ and,   $ \left\{ 
\begin{array}{l}
\lim_{x\rightarrow \pm\infty} u(x)=+\infty, \\
 \lim_{x\rightarrow \pm\infty} u'(x) = \pm\sqrt{\mathrm{e}^{2(c-\varphi(u_0))}-1}. \end{array} \right.$
 \item if $c=+\infty$, 
 $ \lim_{x\rightarrow \pm\Lambda_{u_{0}}} u(x)=+\infty, \quad
 \lim_{x\rightarrow \pm\Lambda_{u_{0}}} u'(x) = \pm\infty.$
\

In particular, if $\Lambda_{u_{0}}<+\infty$, the graph of $u$ is asymptotic to two vertical lines. Moreover,
\item   $\Lambda_{u_{0}}<+\infty$ if and only if $e^{-\varphi}\in L^{1}(]u_0,+\infty[)$,$\left( i.e \int_{u_0}^\infty e^{-\varphi(\tau)} d\tau  < \infty \right)$ .
\item  If   $\Lambda_{\tau}<+\infty$ and $\dot{\varphi}$ is increasing (respectively, decreasing), then $\Lambda_\tau$ is  decreasing (respectively, increasing) in $\tau$.
\end{enumerate}
\end{theorem}

\begin{theorem}[{\rm \cite[Theorem 3.4]{MM}}]\label{t3}
 Let $\varphi:\ ]a,+\infty[\ \longrightarrow \ ]b,c[$,  $ a,b\in \{\R,-\infty\}$,  $c\in \{\R,+\infty\}$ be a strictly decreasing diffeomorphism, then the solution $u$  of \eqref{grafollano}-\eqref{initialcond} is defined in $]-\Lambda_{u_{0}},\Lambda_{u_{0}}[$, $\Lambda_{u_{0}} \in\R^{+}\cup\{+\infty\}$, it is concave,  symmetric about the $y$-axis and has a maximum at $x=0$. Moreover, 
\begin{enumerate}
\item if $c<+\infty$, then 
$\Lambda_{u_{0}}<+\infty $
  and,   $ \left\{ \begin{array}{l} \lim_{x\rightarrow \pm\Lambda_{u_{0}}} u(x)=a, \\
 \lim_{x\rightarrow \pm\Lambda_{u_{0}}} u'(x) = \pm\sqrt{\mathrm{e}^{2(c-\varphi(u_0))}-1}. \end{array} \right.$
\item if $c=+\infty$, then 
$\Lambda_{u_{0}}<+\infty $ if and only if $\int_a^{u_0}\mathrm{e}^{-\varphi(\tau)} d\tau < \infty, $
  and,   $$ \lim_{x\rightarrow \pm\Lambda_{u_{0}}} u(x)=a, \quad
 \lim_{x\rightarrow \pm\Lambda_{u_{0}}} u'(x) = \pm\infty.$$
\end{enumerate}
\end{theorem}

Motivated by their physical interpretation, for each  solution $u$ of \eqref{grafollano}-\eqref{initialcond}  we will refer 
$ {\cal G}^{u_{0}}:= \{(x,y,u(x))\, |\, (x,y)\in I\times \mathbb{R}\}$ as a $[\varphi,\vec{e}_3]$-\textit{catenary cylinder} surface.

\begin{figure}[htb]
\begin{center}
\includegraphics[width=0.3\linewidth]{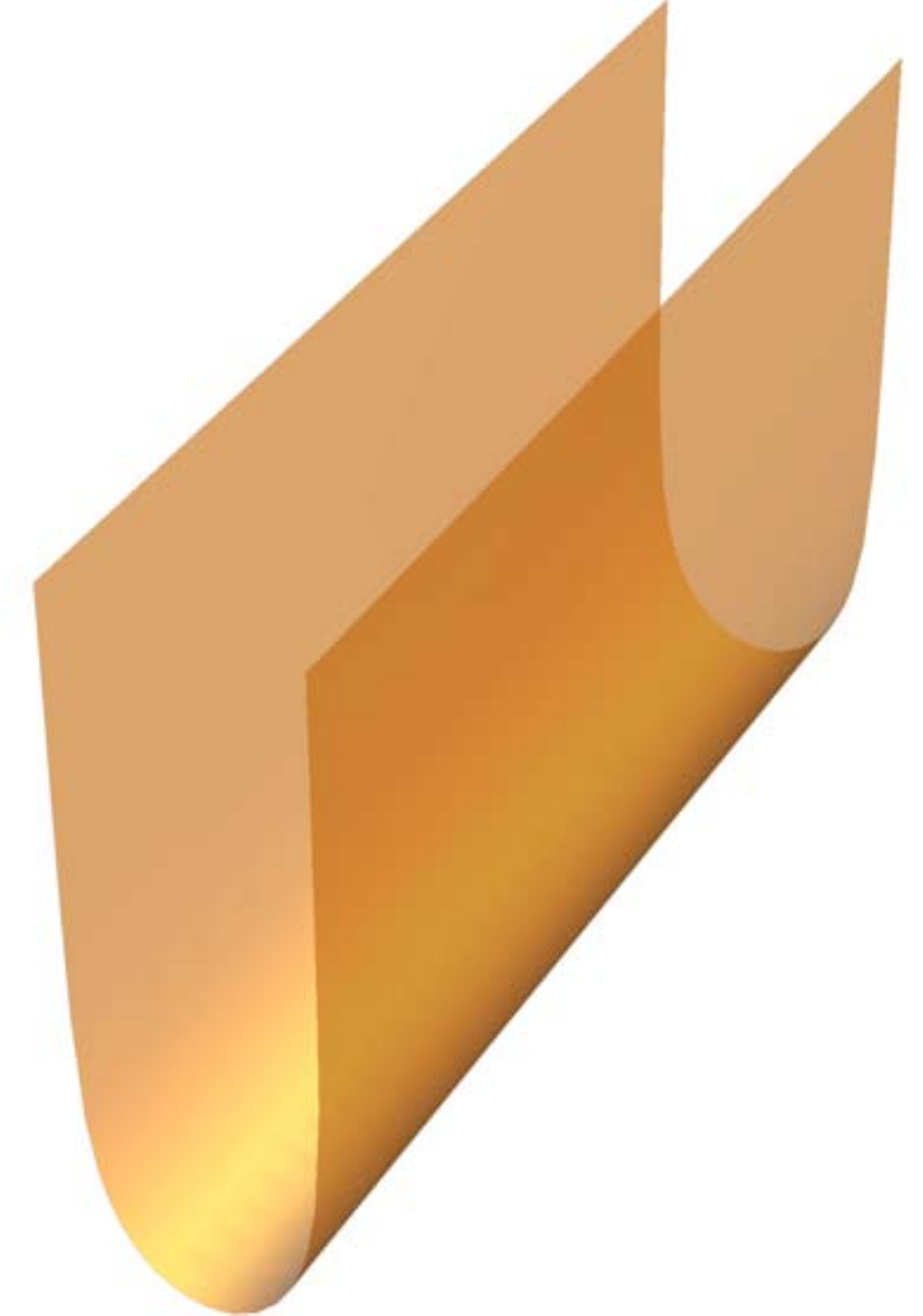}\qquad
\includegraphics[width=0.5\linewidth]{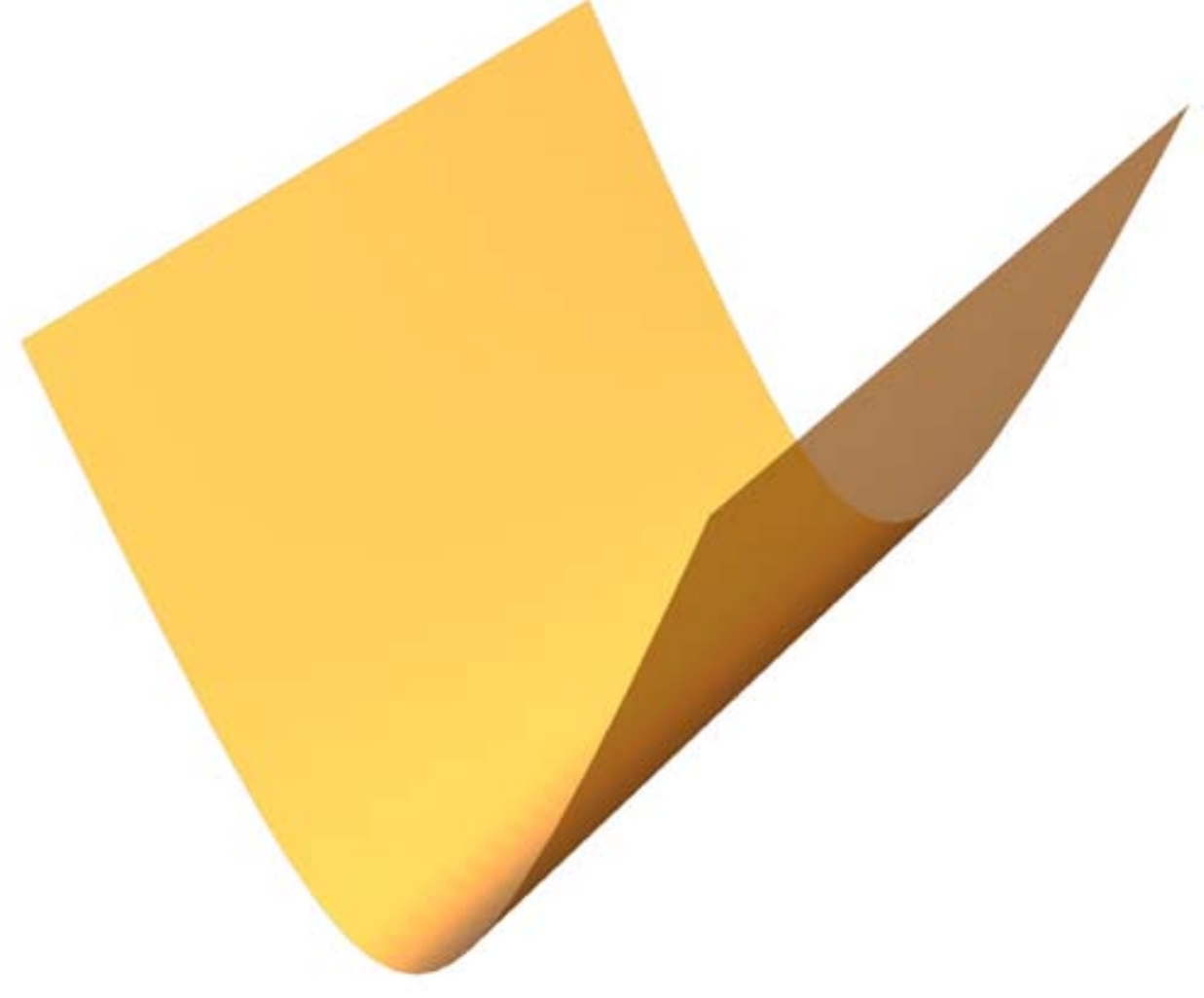}
\end{center}
\caption{$[\varphi,\vec{e}_{3}]$-catenary cylinders with $\dot{\varphi}=1$ and $\dot{\varphi}=1/u^2$, respectively.}
\end{figure}

\subsubsection*{{\sc $\bullet$ Rotationally symmetric solutions}}

In the rotationally symmetric case , the equation \eqref{fe} for $[\varphi,\vec{e}_{3}]$-minimal graphs $u=u(r)$ with $r=\sqrt{x^2+y^2}$ reduces  to the following ODE
\begin{equation}
\label{equrot}
u^{\pp}(r)=(1+u^{\p}(r)^{2})\left(\dot{\varphi}(u)-\frac{u^{\p}(r)}{r} \right),
\end{equation}
where $(^\p)$ denotes the derivative with respect to $r$. Notice that \eqref{equrot} is degenerated and then, the existence and uniqueness of solution around $r=0$ is not guaranteed by standard theory of ordinary differential equations. Moreover, by applying \cite[Theorem 2]{Serrin}, solutions of \eqref{fe} do not have isolated singularities. Consequently, is not a restriction to look for solutions of \eqref{equrot} with the following initial data
\begin{equation}
\label{initialcondrot}
u(0)=u_{0} \ \ , \ \ u^{\p}(0)=0.
\end{equation}
In this case, we can assert (see \cite[Proposition 4.1]{MM}) that  the problem \eqref{equrot}-\eqref{initialcondrot} has a unique solution $u\in C^{2}([0,R])$ for some $R>0$, which depends continuously on the initial data.
\

Now, once the existence of solution is guaranteed, we want to describe  $[\varphi,\vec{e}_3]$-minimal surfaces that  are invariant under  the one-parameter group of rotations that fix the $\vec{e}_{3}$ direction. A such surface with generating curve the arc-lenght parametrized curve $$\gamma(s)=(x(s), 0, z(s)), \qquad   \ \ s\in I\subset\mathbb{R}$$ is given by,
\begin{equation}
\label{param}
\psi(s,t)=\left(x(s)\cos(t), x(s)\sin(t),z(s) \right), \emph{ } (s,t)\in I\times\mathbb{R}.
\end{equation}

The inner normal of $\psi$ writes as
\begin{equation}
\label{Gaussmap}
N(s,t)=\left(-z'(s)\cos(t),-z'(s)\sin(t) ,x'(s) \right),
\end{equation}
and the coefficients of the first and second fundamental form are given by,
\begin{equation}
\label{coefi}
\begin{array}{lll}
&\langle\psi_{s},\psi_{s}\rangle=1,  & \langle\psi_{s},N_{s}\rangle =-\kappa, \\
&\langle\psi_{t},\psi_{t}\rangle=x^{2}, &\langle\psi_{t},N_{t}\rangle = -x \, z', \\ 
&\langle\psi_{s},\psi_{t}\rangle=0, &\langle\psi_{s},N_{t}\rangle =0,
\end{array}
\end{equation}
where $\kappa$ is the curvature of ${\gamma}$ and by $'$ we denote derivative respect to $s$.
Consequently, from \eqref{meancurvature}, \eqref{param}, \eqref{Gaussmap} and \eqref{coefi}, the surface $\psi$ is a $[\varphi,\vec{e}_{3}]$-minimal surface if and only if
\begin{equation}
\label{equz}
\left\{ \begin{array}{l}
x'=\cos(\theta)\\
z'=\sin(\theta),\\
\theta'=\dot{\varphi}(z)\text{cos}(\theta)-\frac{\text{sin}(\theta)}{x},
\end{array}\right.
\end{equation}
where  
$$\theta(s)=\int_0^s \kappa(t) dt.$$
Along this section we will consider that $\varphi:\ ]a,+\infty[\ \rightarrow \R$ is a strictly increasing and convex function, that is
\begin{equation}
\label{conditions}
\dot{\varphi}>0, \quad \ddot{\varphi}\geq0,  \qquad \text{on $]a,+\infty[$}.
\end{equation}
\subsubsection*{{\sc $\bullet$  Globally convex examples:}} Let us  consider the solutions of \eqref{equz} with the following initial conditions,
\begin{equation}
\label{valorinicial} x(0)=0, \qquad z(0)=z_{0}\in ]a,\infty[,  \qquad \theta(0)=0.
\end{equation}
Then, the surface $\psi$ intersects orthogonally the rotation axis and it is globally convex. In fact,  by application of  L'H\^{o}pital's rule we have that $2 \theta'(0)=\dot\varphi(z_0)>0$ and $\gamma$ is a strictly locally convex planar curve around of $s=0$. We assert that   $\theta'(s)>0$ for $s\geq 0$, otherwise from \eqref{valorinicial}, there exists a first value $s_{0}>0$ such that $\theta'(s_{0})=0$ and $\theta''(s_{0})\leq 0$. As $\theta' >0$ on $[0,s_0[$, from \eqref{equz} we have that  $0<2 \theta(s_{0})<\pi$ and by  differentiation  of \eqref{equz}, we get,
$$
\theta''(s_{0})=\frac{\sin(2\theta(s_{0}))}{2}\left(\ddot{\varphi}(z(s_{0}))+\frac{1}{x(s_{0})^2} \right) > 0,
$$
getting to contradiction.
\begin{theorem}[{\rm \cite[Theorem 4.5]{MM}}]
\label{bowl}
Under the conditions \eqref{valorinicial}, the curve ${\gamma}$ is the graph of a strictly convex  symmetric function $u(x)$ defined  on a maximal interval $]-\omega_+,\omega_+[$ which has a minimum at $0$ and $$\lim_{x\rightarrow\pm \omega_+}u(x)=+\infty.$$ 
\end{theorem}
\begin{figure}[htb]
\label{figurebowls}
\begin{center}
\includegraphics[width=0.30\linewidth]{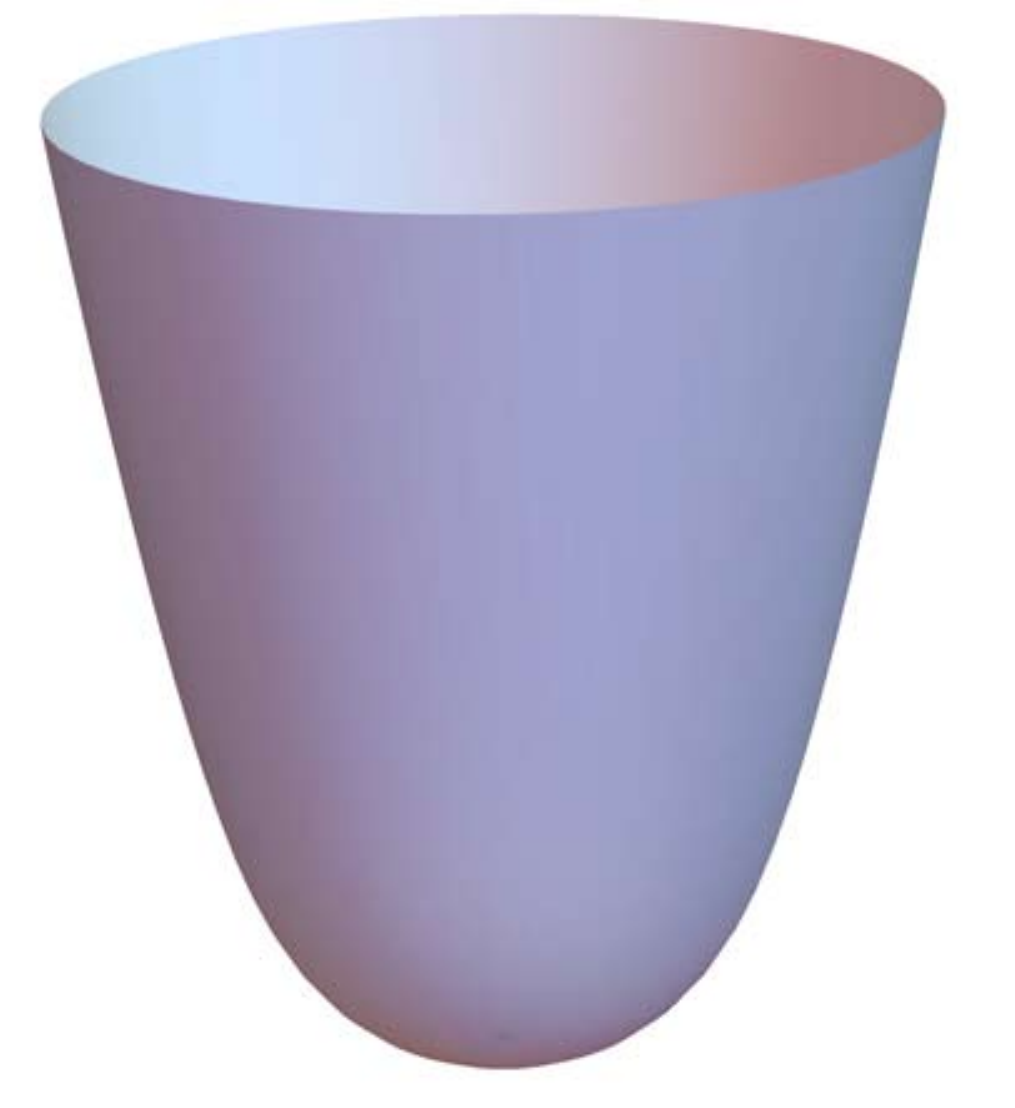}
\qquad \qquad 
\includegraphics[width=0.17\linewidth]{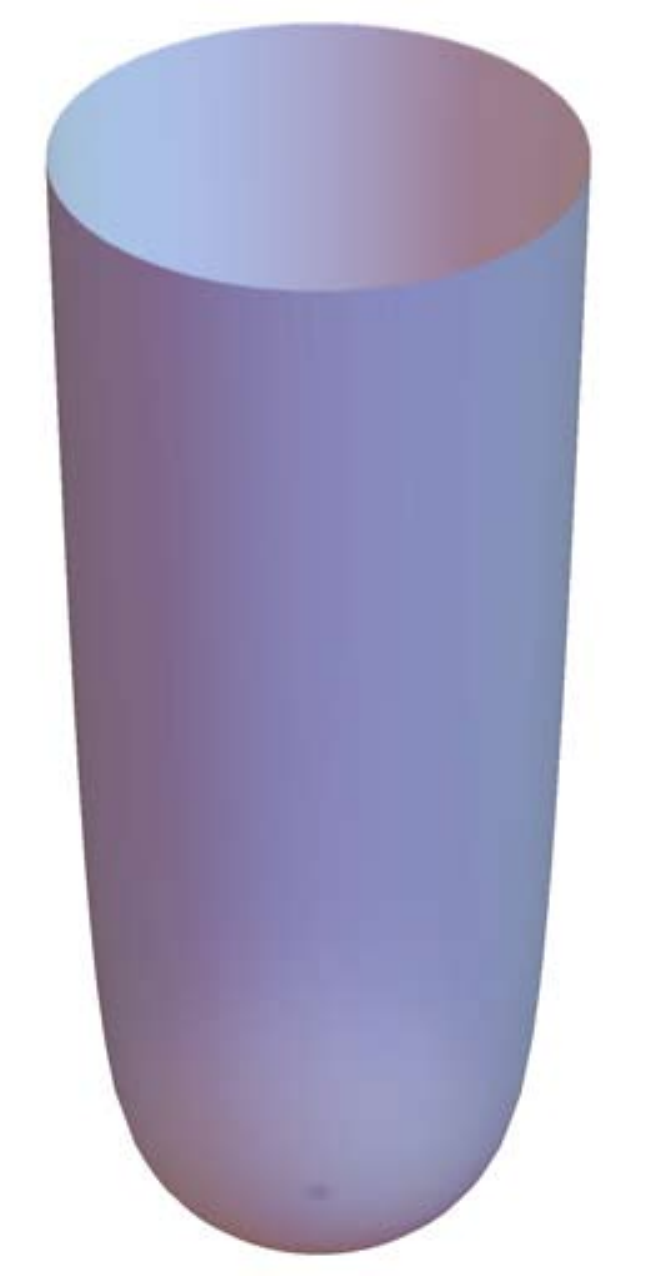}
\end{center}
\caption{$[\varphi,\vec{e}_{3}]$-minimal bowls for $\dot{\varphi}(u)=e^{-1/u}$ and $\dot{\varphi}(u)=u^2$, respectively.}
\end{figure}

 If $\gamma$ is a graph as in Theorem \ref{bowl}, we are going to say that the revolution surface with generating curve $\gamma$ is a \textit{$[\varphi,\vec{e}_{3}]$-minimal bowl} (or a  \textit{$[\varphi,\vec{e}_{3}]$-bowl} in short).

\subsubsection*{{\sc $\bullet$  Non-convex examples:}} Now, we consider the solutions of \eqref{equz} with the following initial conditions,
\begin{equation}
\label{condicionescatenoid} x(0)=x_0>0, \qquad z(0)=z_{0}\in ]a,+\infty[,  \qquad \theta(0)=0.
\end{equation}

From standard theory, the existence and uniqueness of  solution to the  problem \eqref{equz}-\eqref{condicionescatenoid}  is guaranteed. Let $]-s_-,s_+[$ be the maximal interval of existence and consider $\gamma^{+}:=\gamma\big\vert_{[0,s_+[}$ the right branch of $\gamma$. As in Theorem \ref{bowl}, $\gamma^+$ is the graph of a convex function $u=u(x)$ defined on a maximal interval $]x_0,\omega_+[$, such that $$\lim_{x\rightarrow  \omega_+}u(x)=+\infty.$$

 For studying the left branch of $\gamma$ we are going to consider, $\gamma^-(s) = \gamma(-s)$ for $s\in [0,s_-[$. Then, by taking $\overline{x}(s)=x(-s)$,  $\overline{z}(s)=z(-s)$ and  $\overline{\theta}(s)=\theta(-s)+\pi$ for $s\in [0,s_-[$, we have that $\{ \overline{x}, \overline{z},\overline{\theta}\}$ is a solution of \eqref{equz} on $[0,s_-[$ satisfying
$$
\label{condicionescatenoid2} \overline{x}(0)=x_0>0, \qquad \overline{z}(0)=z_{0}\in ]a,+\infty[,  \qquad \overline{\theta}(0)=\pi. $$
Then, we have

\begin{lemma}[{\rm \cite[Section 4.2.2]{MM}}] The following statements hold
\begin{itemize}
\item There exists $s_{0}\in]0,s_{-}[$ such that $2\overline{\theta}(s_{0})=\pi.$
\item If $s\in ]s_{0},s_{-}[$, then $0<2\overline{\theta}(s)<\pi$.
\item $\overline{\theta}$ has a minimum at a point $s_{1}\in ]s_{0},s_{-}[$ and $\overline{\theta}'>0$ on $]s_{1},s_{-}[$.
\item The profile curve $\gamma$ is embedded.
\end{itemize}
\end{lemma}
\begin{figure}[htb]
\begin{center}
\includegraphics[width=0.28\linewidth]{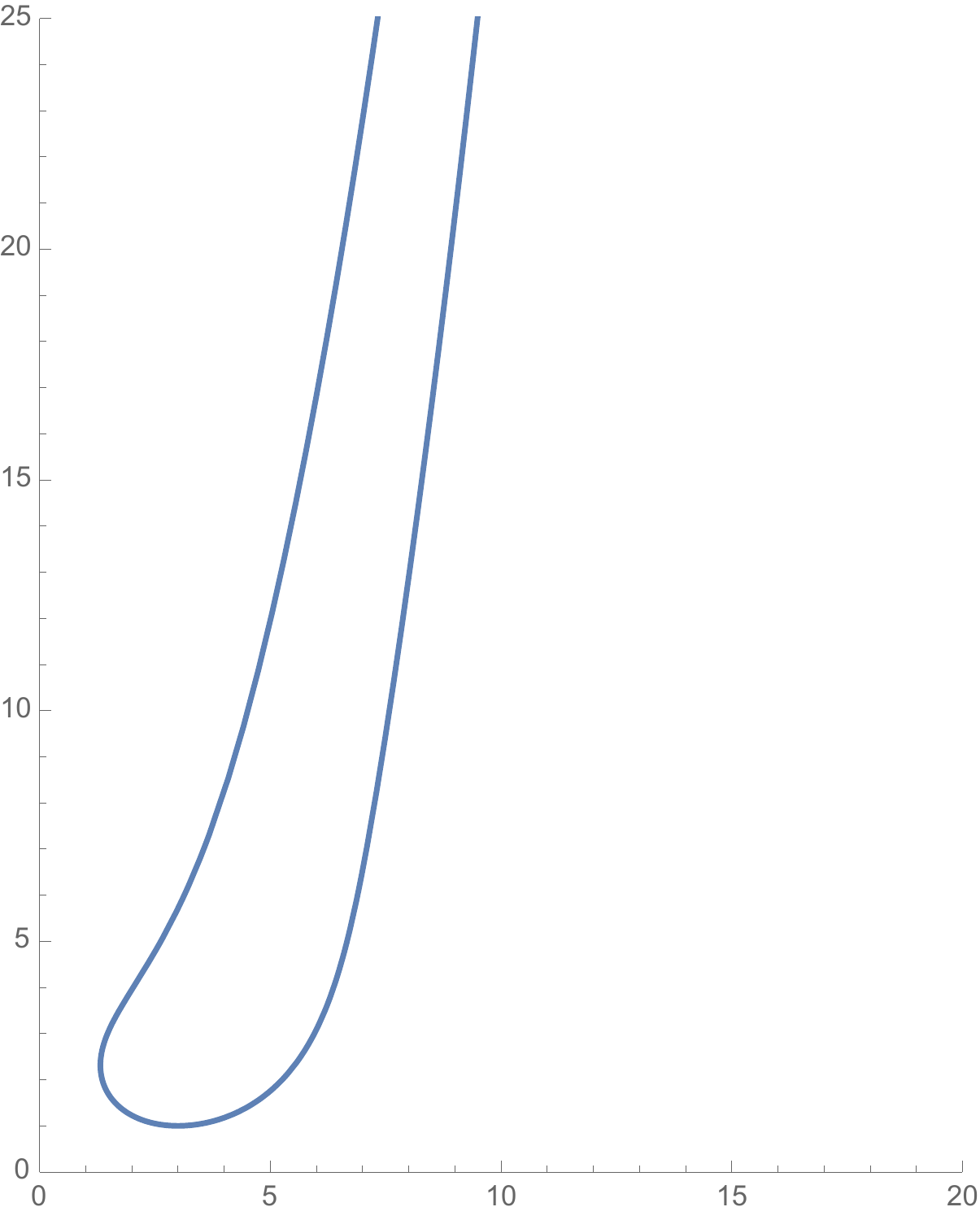}
\quad  
\includegraphics[width=0.30\linewidth]{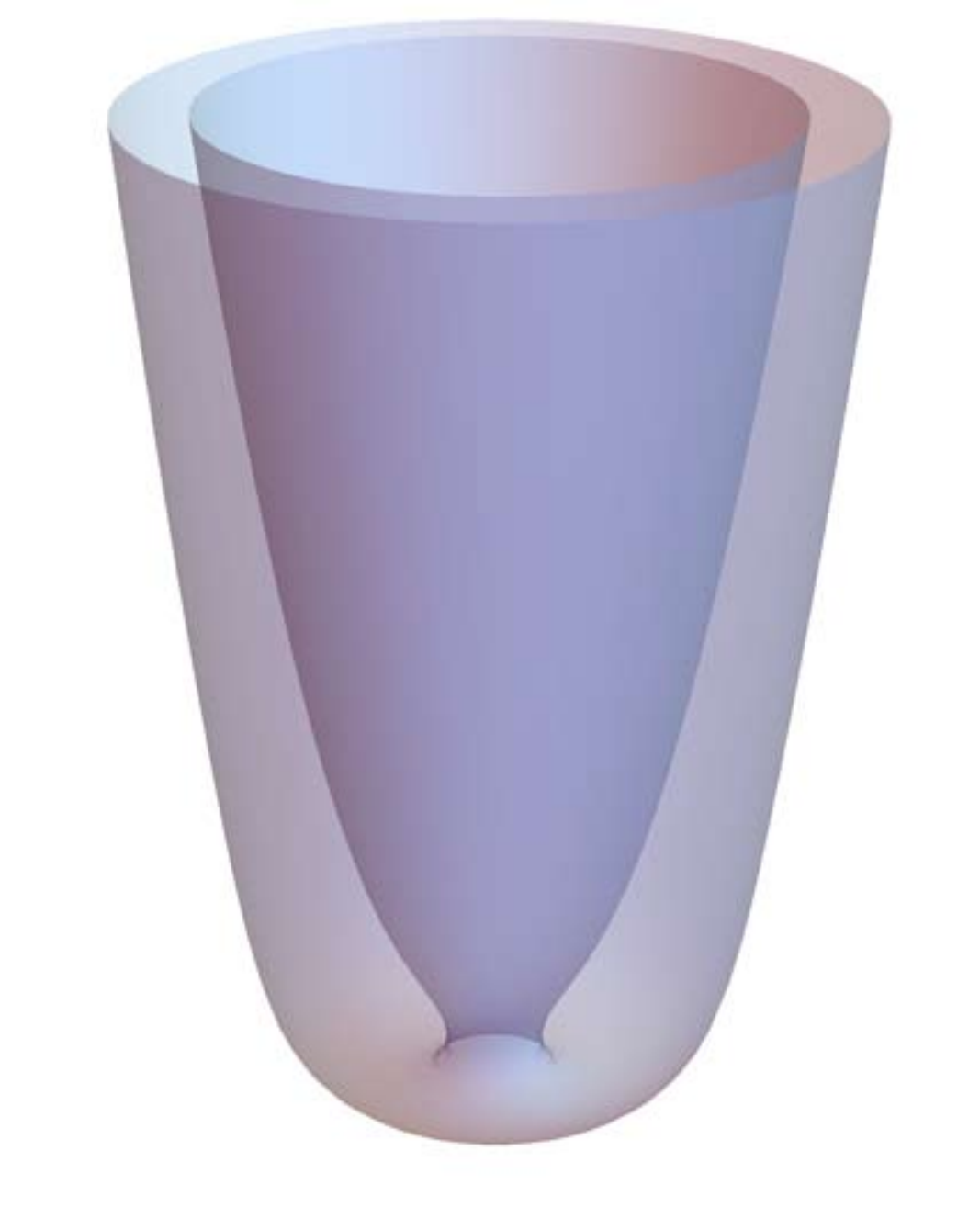}
\end{center}
\caption{$[\varphi,\vec{e}_{3}]$-minimal catenoid with $\dot{\varphi}(u)=e^{-1/u}$.}\label{fcatenoid}
\end{figure}
\begin{theorem}[{\rm \cite[Theorem 4.11]{MM}}]
For every  $x_{0}>0$, there exists a complete embedded rotational $[\varphi,\vec{e}_{3}]$-minimal surface, see Figure \ref{fcatenoid} (right) with the annulus topology whose distance to axis of revolution is $x_{0}$ and whose generating curve $\gamma$ is of winglike type see Figure \ref{fcatenoid} (left). These examples will be called  $[\varphi,\vec{e}_{3}]$-minimal {\sl catenoids}.
\end{theorem}
\subsubsection*{{\sc $\bullet$  The asymptotic behavior}}
One of the questions that we ask ourselves is to know if the rotationally symmetric examples are cylindrically asymptotic or not. Indeed, we may prove
\begin{lemma}[{\rm \cite[Proposition 4.12]{MM}}]Depending on the asymptotic behavior of $\varphi$, we have
\label{cilindro}
\begin{itemize}
\item If $\dot{\varphi}$ has at most a linear growth, then $\omega_{+}=+\infty$ and $\overline{x}_{-}=+\infty$.
\item If $\dot{\varphi}$ growths as $u^{\alpha}$ for some $\alpha>1$, then $\omega_{+},\overline{x}_{-}<+\infty$.
\end{itemize}
\end{lemma}
Motivated by the previous result, we get to control the  asymptotic behavior of the  rotationally symmetric solutions when $\varphi$ has at most quadratic growth. Hence, we are going to assume that $\varphi:]a,+\infty[\rightarrow\mathbb{R}$ is a smooth function satisfying \eqref{conditions} and  with   the following expansion at infinity
\begin{align}\label{series2}
\dot{\varphi}(u) = \Lambda u + \beta + \sum_{n=1}^\infty\frac{a_n}{u^n}, \quad a_n\in \R,
\end{align}
where either $\Lambda>0$ and the first non-vanishing $a_k$ is positive or $\Lambda=0$, $\beta>0$ and the first non-vanishing $a_k$ is negative.

In this case, we can give  explicit formulas for the asymptotic behavior of a rotationally symmetric  solutions of \eqref{equrot} that  generalize  the result of  J. Clutterbuck, O. Schn\"ure and F. Schulze in \cite{CSS} for translating solitons,

\begin{theorem}[{\rm \cite[Theorem A]{MM}}] \label{comportamientoasin}If $\dot{\varphi}$ satisfies \eqref{series2}, then any  rotationally symmetric solution $u$ of  \eqref{fe} has  the following asymptotic behavior,
\begin{itemize}
\item  If  $\Lambda>0$, 
\begin{equation}\label{applineal}
\varphi(u)(r)= C \ e^{\alpha \,r^{2}}  + O(r^2), \quad C>0,
\end{equation}
\item If $\Lambda=0$ and  up to a constant, we have,
\begin{equation}
\label{casoalphacero}
{\cal G}(u)(r)=\frac{r^{2}}{2}-\frac{1}{\beta^2}\log(r)+ {O}(r^{-2}),
\end{equation}
 where ${\cal G}$ is the strictly increasing function given by ${\cal G}(u)=\int_{u_0}^u\frac{d\xi}{\dot{\varphi}(\xi)}$.
 \end{itemize}
\end{theorem}

%
\section{Flat $[\varphi,\vec{e}_{3}]$-minimal surfaces}
Besides the  $[\varphi,\vec{e}_{3}]$-catenary cylinders that we have already described, we mention the following examples of flat $[\varphi,\vec{e}_{3}]$-minimal surfaces: let $\psi:=(x,y,u(x))$, $x\in]-\Lambda_{u_0},\Lambda_{u_0}[$  be a $[\varphi,\vec{e}_3]$-catenary cylinder with $u$ satisfying \eqref{grafollano} and Gauss map,
 $$ N = \frac{1}{\sqrt{1+u'^2}}(u',0,-1). $$
If we rotate the surface by an angle $\theta\in ]0,\pi/2[$ about the $x$-axis and dilate by $1/\cos\theta$, the resulting surface may be written as follows,
$$\widetilde{\psi}=\psi+\frac{1-\cos\theta}{\cos\theta}\langle \psi,\vec{e}_1\rangle \vec{e}_1+(\tan\theta)\vec{e}_1\wedge\psi,$$
where $\vec{e}_1=(1,0,0)$ and 
whose Gauss map  is given by,
\begin{equation}\widetilde{N}=\cos\theta \ N+(1-\cos\theta)\langle N,\vec{e}_1\rangle \vec{e}_1+\sin \theta \ \vec{e}_1\wedge N.\label{ntg}\end{equation}
The mean curvature $\widetilde{H}$ of $\widetilde{\psi}$  verifies
$$\widetilde{H}\ = \cos\theta \ H =\cos \theta \ \dot{\varphi} \langle\vec{e}_3, N\rangle =\dot{\varphi} \langle\vec{e}_3,\widetilde{N}\rangle.$$
Consequently, $\widetilde{\psi}$ is also $[\varphi,\vec{e}_3]$-minimal and we are going to refer these examples as \textit{ tilted $[\varphi,\vec{e}_3]$-catenary cylinders}. 

We have the following classification result,
\begin{figure}[htb]
\begin{center}
\includegraphics[width=0.30\linewidth]{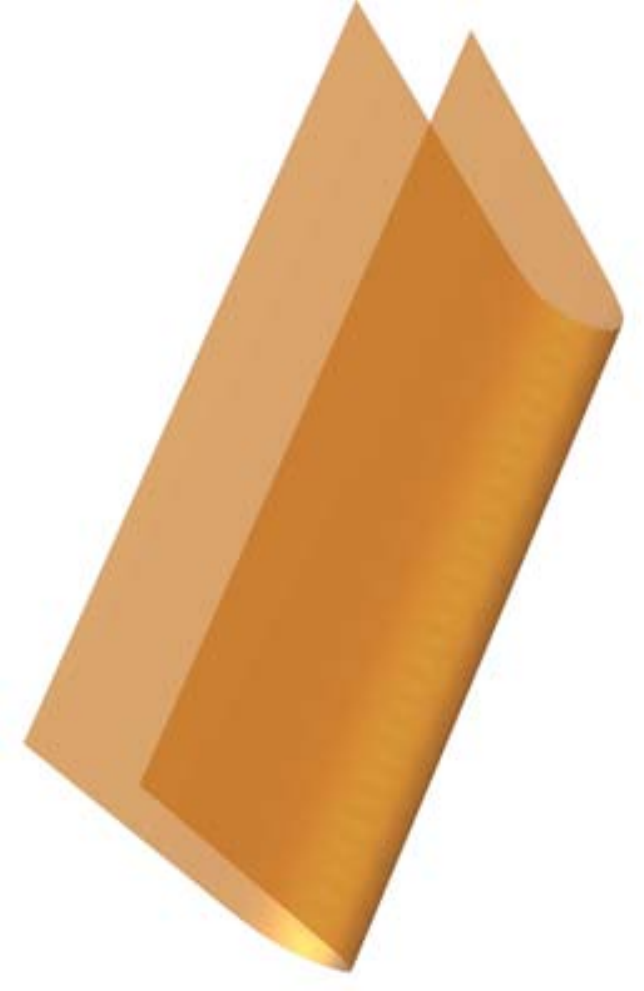}\qquad
\includegraphics[width=0.32\linewidth]{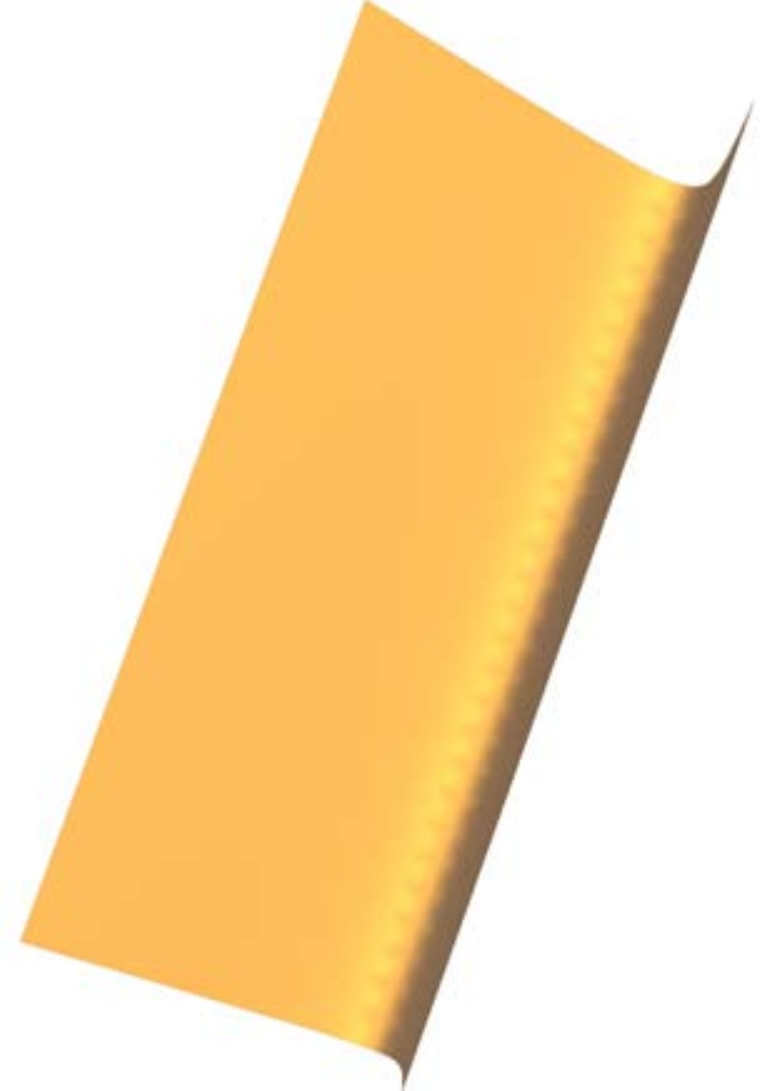}
\end{center}
\caption{Titled $[\varphi,\vec{e}_{3}]$-catenary cylinders with $\dot{\varphi}=1$ and $\dot{\varphi}=1/u^3$, respectively.}
\end{figure}

\begin{theorem}[{\rm \cite[Theorem 3.7]{MM}}]
Let $\Sigma$ be  a complete flat  $[\varphi,\vec{e}_3]$-minimal surface in $\mathbb{R}^{3}$. If $\varphi:\R\rightarrow \R$ is a strictly increasing diffeomorphism, then $\Sigma$ is either a vertical plane or a $[\varphi,\vec{e}_3]$-catenary cylinder (maybe tilted). 
\end{theorem}

\begin{remark}{\rm An analogous clasification was proved for translating solitons by F. Mart\'in, A. Savas-Halilas and K. Smoczyk in \cite{MSHS2} and for singular $\alpha$-minimal surfaces by R. L\'opez in \cite{Rafa1}. }
\end{remark}

%


\section{Mean convex $[\varphi,\vec{e}_{3}]$-minimal surfaces}

In 1983, R. Schoen \cite{Sch} obtained an estimate for the length of the second fundamental form $\mathcal{S}$ of a stable minimal surface $\Sigma$ in a 3-manifold. In particular, in $\mathbb{R}^{3}$, he proved the existence of a constant $C$ such that
$$\vert\mathcal{S}(p)\vert\leq\frac{C}{d_{\Sigma}(p,\partial\Sigma)}, \ \ p\in\Sigma,$$
where $d_{\Sigma}$ stands for  the intrinsic distance of $\Sigma$.
Later,  in 2010, H. Rosenberg, R. Souam and E. Toubiana \cite{RST} obtained an estimate for the length of the second fundamental form, depending on the distance to the boundary, for any stable $H$-surface $\Sigma$ in a complete Riemannian $3$-manifold of bounded sectional curvature $\vert\mathbb{K}\vert\leq\beta < +\infty$. They proved the existence of a constant $C>0$ such that
$$\vert\mathcal{S}(p)\vert\leq\frac{C}{\text{min}\{d_{\Sigma}(p,\partial\Sigma),\pi/2\sqrt{\beta}\}}, \ \ p\in\Sigma.$$
More recently, in 2016, B. White \cite{W1} obtained an estimate for the length of the second fundamental form for minimal surfaces with finite total absolute curvature less than $4\pi$ in $3$-manifolds, depending on the distance to the boundary, on the sectional curvature and on the gradient of the sectional curvature of the ambient space and following C. H. Colding and W. P. Minicozzi  method, \cite{CM,CM1}, J. Spruck and L. Xiao \cite{SX} have also obtained area and curvature bounds for  complete mean convex translating solitons in $\mathbb{R}^{3}$. As application and  using the Omori-Yau maximum principle (see, for example, \cite{AMR}) they have proved  one of the fundamental results in the recent development of translating solitons theory conjectured by X. Wang in \cite{Wang}:
\begin{theorem}{{\rm \cite[Theorem 1.1]{SX}}}
Let $\Sigma \subset \R^3$ be a complete immersed two-sided translating soliton with non-negative mean curvature. Then $\Sigma$ is convex. \label{Thm1.1-SX}
\end{theorem}
In a very recent study, we have extended  the results of \cite{SX} to mean convex $[\varphi,\vec{e}_{3}]$-minimal surfaces. In our case, by mean convex surfaces we will refer to those surfaces with $H \leq 0$ everywhere. More precisely, we will consider mean convex $[\varphi,\vec{e}_{3}]$-minimal oriented surfaces $\Sigma$ with empty boundary in $\R^3_\alpha=\{p\in\R^3\ \vert \ \langle p,\vec{e}_{3}\rangle>\alpha\}$, where $\varphi: \R \longrightarrow \R$ is a smooth function  satisfying
\begin{equation} \dot{\varphi} >0, \quad \ddot{\varphi}\geq 0 \quad \text{on $]\alpha,+\infty[$}. \label{c1}
\end{equation}
In order to describe our results, let  us consider $F_{t}$  the normal variation associated to a compactly supported variational vector field on the normal bundle of $\Sigma$.  Then, see \cite[Appendix]{CMZ}, the second variation of $\mathcal{A}^{\varphi}$ is given by,
$$\frac{d^{2}}{dt^{2}}\bigg\vert_{t=0}\mathcal{A}^{\varphi}(F_{t}(\Sigma))=-\int_{\Sigma}\varrho\mathcal{L}_{\varphi}(\varrho)\, d\Sigma \ \ \text{ for any } \varrho\in C^{\infty}_{0}(\Sigma),$$
where $\mathcal{L}_{\varphi}$ is a gradient Schr\"odinger's type operator defined on $C^{2}(\Sigma)$ and given by
$$\mathcal{L}_{\varphi}(\cdot)=\Delta^{\varphi}(\cdot)+(\vert\mathcal{S}\vert^{2}-\ddot{\varphi}\eta_{3}^{2})(\cdot)$$
where $\Delta^\varphi (\cdot ) = \Delta  (\cdot ) + \langle \nabla \varphi, \nabla (\cdot )  \rangle$  and $\eta_3=\langle N,\vec{e}_3\rangle$. 

 As usual, we will say that $\Sigma$ is stable if and only if for any compactly supported smooth function $\varrho$, it holds that
\begin{equation}
\label{stable}
-\int_{\Sigma}\varrho\mathcal{L}_{\varphi}(\varrho)e^{\varphi}\, d\Sigma\geq 0.
\end{equation}
It is interesting to mention the following fact:
\begin{proposition}[{\rm \cite[Proposition 4.4]{MMJ}}]\label{stableHpositiva}
Let $\varphi:]\alpha,+\infty[\  \rightarrow \R$ be a regular function satisfying \eqref{c1}  and $\Sigma$ be an oriented $[\varphi,\vec{e}_{3}]$-minimal immersion in $\mathbb{R}^{3}$ with $H\leq 0$. Then, $\Sigma$ is stable.
\end{proposition}

\begin{remark}{\rm 
The existence of stable surfaces is not guaranteed for any function $\varphi$. X. Cheng, T. Mejia and D. Zhou \cite{CMZ} proved that if the Ilmanen's space is complete and $\ddot{\varphi}\leq-\epsilon<0$, for some $\epsilon>0$, then there are not stable surfaces without boundary and with finite weighted area.}
\end{remark}
From Proposition \ref{stableHpositiva} and following the same method as in \cite{SX}, we get the following area stimate
\begin{theorem}[{\rm \cite[Proposition 4.10]{MMJ}}]
\label{boundnessarea}
Let  $\Sigma$ be a  $[\varphi,\vec{e}_{3}]$- immersion in $\mathbb{R}^{3}_\alpha$ with $H\leq 0$ and $\varphi$ satisfying \eqref{c1} and \begin{equation}
\Gamma :=\sup_{]\alpha,+\infty[}  (2 \ddot{\varphi} - \dot{\varphi}^2 )< +\infty. \label{c2}.
\end{equation} If  $2\rho\,\dot{\varphi}(\rho+\mu(p))<\log(2)$ and $ \sqrt{\vert \Gamma \vert} \ \rho < 1$, then the geodesic disk $\mathcal{D}_{\rho}(p)$ of radius $\rho$  centered at $p$ is disjoint from the cut locus of $p$ and 
\begin{equation}
\label{cotarea}
\mathcal{A}(\mathcal{D}_{\rho}(p))< 4\pi\rho^{2},
\end{equation}
 where $\mathcal{A}(\cdot)$ is the intrinsic area of $\Sigma$ in $\R^3$.
\end{theorem}
But, for obtaining curvature bounds we need a better control at infinity of the function $\varphi$ . To be more precise, we are going to consider that   $ z \mapsto \frac{\dot{\varphi}(z)}{z}$ is analytic at $+\infty$; i.e.,
$\dot{\varphi}$ has the following series expansion at $+\infty$:
\begin{align}
& \dot{\varphi}(u)= \Lambda \, u+\beta+\sum_{i=1}^{\infty}\frac{c_{i}}{u^{i}}, \quad \text{$u$ large enough},\label{cc3}
\end{align}
with $\Lambda\geq 0$ and $\beta>0$ if $\Lambda=0$.
\begin{remark}{\rm 
It is worth to note that  condition \eqref{cc3} implies \eqref{c2}. Besides of a natural extension of the best known examples,  conditions \eqref{c1} and \eqref{cc3}  are interesting because under these assumptions it is possible to know explicitly the asymptotic behavior of rotational and translational invariant examples (see \S 2). }
\end{remark}
Bearing in mind Proposition \ref{stableHpositiva}, Theorem \ref{boundnessarea} and the compactness Theorem  2.1 of B. White \cite{W2} for minimal surfaces together with the clasification of all complete translating soliton graphs in $\mathbb{R}^{3}$ (see \cite{HIMW}) we can prove the following Blow-up result:
\begin{theorem}[{\rm \cite[Theorem 4.13]{MMJ}}]
\label{blowup}
Let   $\Sigma$ be a properly embedded $[\varphi,\vec{e}_{3}]$-minimal surface in $\mathbb{R}^{3}_\alpha$ with $H\leq 0$,  locally bounded genus and $\varphi$ satisfying \eqref{c1} and \eqref{cc3}. Consider any sequence $\{\lambda_{n}\}\rightarrow+\infty$ and suppose that there exists a sequence $\{p_{n}\}$ in $\Sigma$ such that  $\{\dot{\varphi}(\mu(p_{n}))/\lambda_{n}\}\rightarrow C$ for some constant $C\geq 0$ .  Then, after passing to a subsequence,   $\Sigma_n=\lambda_n ( \Sigma-p_n)$  converge smoothly to
\begin{enumerate}
\item a plane  when $C=0$,
\item one of the following translating solitons when $C>0$:
\begin{enumerate}
\item[(a)] vertical plane,
\item[(b)] grim reaper surface,
\item[(c)] titled grim reaper surface,
\item[(d)] bowl soliton,
\item[(e)] $\Delta$-Wing translating soliton.
\end{enumerate}
\end{enumerate}
\end{theorem}
 From Theorem \ref{blowup} and  by combining the methods of H. Rosenberg, R. Souam and E. Toubiana \cite{RST} and J. Spruck and L. Xiao \cite{SX}, we have the following  curvature estimates, 
\begin{theorem}[{\rm \cite[Theorem A]{MMJ}}]\label{mt1} Let $\Sigma$ be a  properly embedded $[\varphi,\vec{e}_{3}]$-minimal surface in $\R^3_\alpha$  with non-positive mean curvature, locally bounded genus and  $\varphi: \R \rightarrow \R$  satisfying \eqref{c1} and \eqref{cc3}. Then $\vert\mathcal{S}\vert/\dot{\varphi}$ is bounded on  $\Sigma$. In particular,  if  $\Lambda=0$,  $\vert\mathcal{S}\vert$ is  bounded and if $\Lambda\neq 0$, $\vert\mathcal{S}\vert$ may go to infinity but with at most a linear growth in height. 
\end{theorem}
 As we have already mentioned, J. Spruck and  L. Xiao \cite{SX} proved that any complete translating soliton in $\mathbb{R}^{3}$ with $H\leq 0$ is convex. Later,  D. Hoffman, T. Ilmanen,  F. Mart\'in and B. White in \cite{HIMW2}, by using the Omori-Yau's maximum principle for the Laplacian $\Delta$ \cite[Theorem 3.2]{AMR}, have obtained a more simplified proof of it. 

As an extension of this result, we have
%
\begin{theorem}[{\rm \cite[Theorem B]{MMJ}}]\label{mt2} Let $\Sigma$ be a  properly embedded $[\varphi,\vec{e}_{3}]$-minimal surface in $\R^3_a$  with non positive mean curvature, locally bounded genus and  $\varphi: \R \rightarrow \R$  satisfying \eqref{c1},  \eqref{cc3} and $\dddot{\varphi}\leq 0$ on $]a,+\infty[$. Then $\Sigma$ is convex if and only if the function $\Lambda K$ is bounded from below, where $K$ is the Gauss curvature.
\end{theorem}
\begin{remark} {\rm The main tool in the proof of Theorem \ref{mt2} is the Omori-Yau's maximum principle for the drift Laplacian $\Delta^{\varphi}$ \cite[Theorem 3.2]{AMR} and  it is remarkable that the condition \eqref{cc3} on $\varphi$ is essential for proving that this maximum principle can be applied (see \cite[Theorem 5.1]{MMJ}).
}
\end{remark}
%

\section{Uniqueness of  Dirichlet's  problems at infinity}
Despite of this large family of surfaces and the very  general conditions on $\varphi$, it is possible to prove  uniqueness of bowls and $[\varphi,\vec{e}_{3}]$-catenary cylinders from their asymptotic behavior.
\subsubsection*{{\sc $\bullet$ Uniqueness of  $[\varphi,\vec{e}_{3}]$-bowls}}
Let $\varphi:]a,+\infty[\rightarrow\mathbb{R}$ , $a\in\mathbb{R}\cup\{-\infty\}$ be a strictly increasing convex smooth function satisfying \eqref{series2} and let $\Sigma$ be a complete connected properly embedded $[\varphi,\vec{e}_{3}]$-minimal surface  and $\mathfrak{D}^3=\mathbb{R}^{2}\times]a,+\infty[$. From Theorem \ref{comportamientoasin}, it is natural to say  that an end of  $\Sigma$  is {\sl smoothly asymptotic} to a  rotational-type example if  $\Sigma$ can be expressed outside a Euclidean ball as a vertical graph of  a function $u_\Sigma$ so that,
\begin{equation}\label{defasym}
\varphi (u_{\Sigma})(x)= C \,e^{\Lambda \,|x|^2}  + O\left(|x|^2\right),\quad \text{if}\quad \Lambda>0,
\end{equation}
where $C$ is a positive constant  and, up to a constant, 
\begin{equation}\label{defasym2}
{\cal G}(u_{\Sigma})(x)= \frac{|x|^2}{2}-\frac{1}{\beta^2}\log(|x|)+ {O}\left(|x|^{-2}\right), \text{ if } \Lambda =0 \text{ and } \beta>0.
\end{equation}
where ${\cal G}$ is the strictly increasing function given by ${\cal G}(u)=\int_{u_0}^u\frac{d\xi}{\dot{\varphi}(\xi)}$.

\

Under these conditions the following result holds, 
\begin{theorem}[{\rm \cite[Theorem B]{MM}}]
\label{unicidad} Let $\varphi:]a,+\infty[\rightarrow\mathbb{R}$ , $a\in\mathbb{R}\cup\{-\infty\}$ be a strictly increasing convex smooth function satisfying \eqref{series2} and 
 $\Sigma$ be a complete properly embedded  $[\varphi,\vec{e}_3]$-minimal surface in $\mathbb{R}^{3}$ with a single end that is smoothly asymptotic to a  $[\varphi,\vec{e}_3]$-minimal bowl. Then the surface $\Sigma$ is a  $[\varphi,\vec{e}_3]$-minimal bowl.

\end{theorem}
\begin{remark}{\rm 
The proof of  Theorem \ref{unicidad} is based on the use of the Alexandrov reflection principle (see \cite{Al}) to prove that the surface is symmetric with respect to any vertical plane  through the origin. Although this principle is applied in a  standard way, it is crucial in the proof (see \cite[Lemmas 6.3 and 6.4]{MM}) to show that it is possible to start the reflexion  respect to any vertical plane far enough from the origin.}
\end{remark}
\begin{remark} {\rm In the particular case of translating solitons, Theorem \ref{unicidad} was proved in \cite[Theorem A]{MSHS2}.}
\end{remark}

\subsubsection*{{\sc $\bullet$ Uniqueness of  $[\varphi,\vec{e}_{3}]$-catenary cylinders:}}

Let $\varphi: ]a,+\infty[\rightarrow]b,c[$ $a,b\in \mathbb{R}\cup\{-\infty\}$ , $c\in\mathbb{R}\cup\{+\infty\}$ be a strictly increasing diffeomorphims such that $e^{-\varphi} \in L^{1}(]a,+\infty[)$ and let us consider $\mathcal{G}^{u_0}$, $u_0\in ]a,+\infty[$  one of the  $[\varphi,\vec{e}_{3}]$-catenary cylinders described in $\S 2$.

If $\Sigma$ is a complete connected $[\varphi,\vec{e}_{3}]$-minimal graph and $\Pi(0)$ is the vertical plane through the origin orthogonal to $\vec{e}_{1}$, we will say that a smooth surface $\Sigma$ is $\mathcal{C}^{k}$-\textit{asymptotic} to the \textit{right part} 
$$\mathcal{G}^{u_0}_{+}(0)=\{p\in \Sigma \,:\,\langle p,\vec{e}_1\rangle \, \geq 0\}$$ 
of $\mathcal{G}^{u_0}$  if for any $\varepsilon>0$ there exists $\delta>0$ such that $\Sigma$  can parametrized as a graph over $\mathcal{G}^{u_0}$ as follows 
\begin{equation}
\label{p1}
\widetilde{F}:T_{\delta,u_0}^{+}\subset\mathbb{R}^{2}\rightarrow\mathbb{R}^{3} \ \  \ \ \widetilde{F}=F+\overline{u}\, N_{F},
\end{equation}
where $T_{\delta,u_0}^{+}:=]\Lambda_{u_{0}}-\delta,\Lambda_{u_{0}}[\times\mathbb{R}$, $F(x_{1},x_{2})=(x_{1},x_{2},u(x_{1}))$ parametrizes $\mathcal{G}^{u_0}$ on $T_{\delta,u_0}^{+}$, $u$ is a solution of \eqref{grafollano} with $u(0)=u_0$, $\overline{u}:T_{\delta,u_0}^{+}\rightarrow\mathbb{R}$ is a function in $\mathcal{C}^{k}(T_{\delta,u_0}(+))$ such that
\begin{equation}
\label{p2}
\sup_{T_{\delta,u_0}^{+}}\vert\overline{u}\vert <\varepsilon \ \ , \ \ \sup_{T_{\delta,u_0}^{+}}\vert D^{j}\overline{u}\vert<\varepsilon, \text{ for any } j\in\{1,\cdots,k\}.
\end{equation}
and $N_{F}$ is the downwards unit normal of $\mathcal{G}^{u_0}$. 

Analogously,  we will say that a smooth surface $\Sigma$ is $\mathcal{C}^{k}$-\textit{asymptotic} to left \textit{left part} 
$$\mathcal{G}^{u_0}_{-}(0)=\{p\in \Sigma \,:\,\langle p,\vec{e}_1\rangle \, \leq 0\}$$ 
of $\mathcal{G}^{u_0}$  if for any $\varepsilon>0$ there exists $\delta>0$ such that $\Sigma$  can parametrized as a graph over $\mathcal{G}^{u_0}$ as follows \begin{equation}
\label{p3}
\widetilde{F}:T_{\delta,u_0}^{-}\subset\mathbb{R}^{2}\rightarrow\mathbb{R}^{3} \ \  \ \ \widetilde{F}=F+\overline{u}\, N_{F},
\end{equation}
where $T_{\delta,u_0}^{-}:=]-\Lambda_{u_{0}},-\Lambda_{u_{0}}+\delta[\times\mathbb{R}$, $F(x_{1},x_{2})=(x_{1},x_{2},u(x_{1}))$ parametrizes $\mathcal{G}^{u_0}$ on $T_{\delta,u_0}^{-}$, $u$ is a solution of \eqref{grafollano} with $u(0)=u_0$, $\overline{u}:T_{\delta,u_0}^{-}\rightarrow\mathbb{R}$ is a function in $\mathcal{C}^{k}(T_{\delta,u_0}^-)$ such that
\begin{equation}
\label{p4}
\sup_{T_{\delta,u_0}^{-}}\vert\overline{u}\vert <\varepsilon \ \ , \ \ \sup_{T_{\delta,u_0}^{-}}\vert D^{j}\overline{u}\vert<\varepsilon, \text{ for any } j\in\{1,\cdots,k\}.
\end{equation}

In particular, we say that $\Sigma$ is $\mathcal{C}^{k}$-\textit{asymptotic} to $\mathcal{G}^{u_0}$ if and only if $\Sigma$ is $\mathcal{C}^{k}$-asymptotic to the both branches $\mathcal{G}^{u_0}_{+}(0)$ and $\mathcal{G}^{u_0}_{-}(0)$. Moreover, a smooth surface $\Sigma$ is called $\mathcal{C}^{k}$-\textit{asymptotic to $\mathcal{G}^{u_0}$, outside a cylinder}, if  there exists a solid cylinder $\mathfrak{c}$ whose axis is $\mathcal{G}^{u_0}\cap\Pi(0)$ and the set $\Sigma-\mathfrak{c}$ consists of two connected components $\Sigma_{1}$ and $\Sigma_{2}$ which are $\mathcal{C}^{k}$-asymptotic to $\mathcal {G}^{u_0}_{+}(0)$ and $\mathcal{G}^{u_0}_{-}(0)$, respectively.

\

As a consequence of the compactness result \cite[Theorem 3.4]{MJ} and \cite[Lemma 4.3]{MJ}, we have
\begin{proposition}[{\rm \cite[Proposition 4.5]{MJ}}]
\label{comportamientocurvas}
Let $\varphi:]a,+\infty[\rightarrow ]b,c[$ , $a,b\in\mathbb{R}\cup\{-\infty\}$ , $c\in\mathbb{R}\cup\{+\infty\}$ be a convex strictly increasing diffeomorphism with $e^{-\varphi}\in L^{1} (]a,+\infty[)$ and $\Sigma$ be a connected $[\varphi,\vec{e}_{3}]$-minimal immersion $\mathcal{C}^{\infty}$-asymptotic to $[\varphi,\vec{e}_{3}]$-catenary cylinder $\mathcal{G}^{h}$, outside a cylinder, for some $h\in ]a,+\infty[$. For any sequence of points $\{(p_{1,n},p_{2,n},p_{3,n})\}$ of $\Sigma$ such that $\{p_{2,n}\}$ diverges and $\{p_{3,n}\}$ is bounded, the sequence  $\{\Sigma_{n}=\Sigma-(0,p_{2,n},0)\}_{n\in\mathbb{N}}$ converge smoothly, after subsequence, to some $[\varphi,\vec{e}_{3}]$-catenary cylinder with the same asymptotic behaviour that $\mathcal{G}^{h}$.
\end{proposition}
From Proposition \ref{comportamientocurvas},  there exists  $n_{0}$ large enough such that for any $n\geq n_{0}$, each $\Sigma_{n}$ can be parametrized as a graph over some $\mathcal{G}^{u_0'}$ in $\mathcal{T}_{\delta,u_0,n}^{+}$ (respectively $\mathcal{T}_{\delta,u_0,n}^{-}$) by
\begin{equation}
\label{p5}
\widetilde{F}_{n}:\mathcal{T}_{\delta,u_0,n}^{+}\rightarrow\mathbb{R}^{3} \ \ \ \ \widetilde{F}_{n}=F+\overline{u_{n}}N_{F},
\end{equation}
where  $\mathcal{T}_{\delta,u_0,n}^{+}=]-\Lambda_{u_0'}+\delta,\Lambda_{u_0'}-\delta[\times ]m_{1,n},m_{2,n}[\rightarrow\mathbb{R}$ (respectively,$\mathcal{T}_{\delta,u_0,n}^{-}=]-\Lambda_{u_0'}+\delta,\Lambda_{u_0'}-\delta[\times ]-m_{2,n},-m_{1,n}[$) with $\Lambda_{u_0'}=\Lambda_{u_{0}}$ due to the asymptotic behavior, $\delta>0$ only depends of $n$, $\{m_{1,n}\}_{n\in\mathbb{N}},\{m_{2,n}\}_{n\in\mathbb{N}}$ are strictly monotonous sequences with $m_{1,n}<m_{2,n}$ and each $\overline{u_{n}}:\mathcal{T}_{\delta,u_0,n}^{+}\rightarrow\mathbb{R}$ (respectively,  $\overline{u_{n}}:\mathcal{T}_{\delta,u_0,n}^{-}\rightarrow\mathbb{R}$) is a smooth function satisfying the following inequalities,
\begin{equation}
\label{p6}
\text{sup}_{\mathcal{T}_{\delta,u_0,n}^{+}}\vert\overline{u_{n}}\vert<\varepsilon, \ \ \text{sup}_{\mathcal{T}_{\delta,u_0,n}^{+}}\vert D^{j}\overline{u_{n}}\vert<\varepsilon, \ \ \text{for any } j\in\mathbb{N}.
\end{equation}
By using the above inequalities it is possible to prove that  $$\frac{\eta_2}{\eta_3}:=\frac{<N,\vec{e}_2>}{<N,\vec{e}_2>},$$  goes to zero at infinity an then, there exists an interior point where the function $\eta_{2}/\eta_{3}$ attains either a local minimum in $\{p\in\Sigma: \eta_{2}(p)<0\}$ or a local maximum in $\{p\in\Sigma: \eta_{2}(p)>0\}$. But then, it is possible to deduce that $\eta_{2}$ vanishes everywhere and so $\Sigma$ is invariant under translations in the direction $\vec{e}_{2}$ which gives

\begin{theorem}[{\rm \cite[Theorem 1.4]{MJ}}]
\label{uni}
Let $\varphi:]a,+\infty[\rightarrow ]b,c[$, $a,b\in\mathbb{R}\cup\{-\infty\}$ and $c\in\mathbb{R}\cup\{+\infty\}$ be a strictly increasing convex diffeomorphism such that $e^{-\varphi}\in L^{1}(]a,+\infty[)$ and bounded quotient $\ddot{\varphi}/\dot{\varphi}$. If $\Sigma$ is a complete connected $[\varphi,\vec{e}_{3}]$-minimal graph $\mathcal{C}^{\infty}$-asymptotic to $[\varphi,\vec{e}_{3}]$-catenary cylinder $\mathcal{G}^{u_0}$, outside a cylinder, for some $h\in ]a,+\infty[$, then $\Sigma$ coincides with some $[\varphi,\vec{e}_{3}]$-catenary cylinder with the same behaviour that $\mathcal{G}^{u_0}$.
\end{theorem}
%

\begin{remark}{\rm 
F. Mart\'in, J. P\'erez-Garc\'ia, A. Savas-Halilaj and K. Smoczyk \cite{MSHS1} proved that, if $\Sigma$ is a properly embedded translating soliton with locally bounded genus and $C^{\infty}$-topology to two vertical planes outside a cylinder, then $\Sigma$ must be coincide with some grim reaper translating soliton.}
\end{remark}
\section{Weierstrass' type representation}
Let $\varphi: I\rightarrow \mathbb{R}$ be a smooth function on a real open interval $I\subseteq \mathbb{R}$ 
 and  $ \psi:\Sigma\rightarrow\mathfrak{D}^3=\R^2\times I\subseteq\R^3$ be an immersion. Consider a local conformal parameter $\zeta=u+iv$ of $\Sigma$ on an open simply connected domain ${\cal U} \subset \mathbb{C}$ such that the induced metric $ds^2$  writes
\begin{equation} ds^2 := \lambda^2 (du^2 +dv^2) = \lambda^2 |d\zeta|^2\label{imetric}
\end{equation} and set, as usual, the  Wirtinger's operators  by, 
$$
\partial_{\zeta}=\frac{1}{2}\left(\partial_{u}-i\partial_{v}\right), \quad \partial_{\overline{\zeta}}=\frac{1}{2}\left(\partial_{u}+i\partial_{v}\right).
$$ 
If $\psi_\zeta=\mathrm{e}^{\frac{1}{2} \varphi}(f,g,h)$ then, the conformality conditions write
\begin{align}
& \lambda^2 = 2\vert\psi_\zeta\vert^2=2 \, \mathrm{e}^\varphi\left(\vert f\vert^{2}+\vert g\vert^{2}+\vert h\vert^{2}\right)\label{con1}\\
&f^{2}+g^{2} + h^{2}=0,\label{con2}
\end{align}
and by an straightforward computation we obtain that $\psi$ is a $[\varphi, \vec{e}_3]$-minimal immersion if and only if 
\begin{equation}
\label{sistema2}
\mathrm{e}^{\frac{1}{2} \varphi}h_{\overline{\zeta}}=\frac{1}{2}\dot{\varphi}\left(\vert f\vert^{2} +\vert g\vert^{2}\right), \quad \mathrm{e}^{\frac{1}{2} \varphi}f_{\overline{\zeta}}=-\frac{1}{2}\dot{\varphi}\overline{f}h,  \quad
\mathrm{e}^{\frac{1}{2}\varphi} g_{\overline{\zeta}}=-\frac{1}{2}\dot{\varphi}\overline{g}h
\end{equation}
Now, if we introduce the complex functions 
\begin{equation}
\label{complexfunctions}
F=f-i\ g \quad \text{ and } \quad G=\frac{h}{F}, 
\end{equation}
from \eqref{con1} and \eqref{con2}, we have that  $G$ is a smooth map into the Riemann sphere and if $G$ is not constant,
\begin{equation}
\label{coordiso}
h=FG \ \ , \ \ f=\frac{1}{2}F(1-G^{2}) \ \ , \ \ g=\frac{i}{2}F(1+G^{2}).
\end{equation}
Moreover, the Gauss map $N$ of $\psi$ in the Euclidean space $\mathbb{R}^3$ is given in terms of $G$ as
\begin{align*}
 N:= \left(\frac{2G}{1+\vert G\vert^{2}}, \frac{1-\vert G\vert^{2}}{1+\vert G\vert^{2}} \right).
\end{align*}
 We are going to say also  that $G$ is the \textit{Gauss map} of  $\psi$.  

\

By using \eqref{con1} and \eqref{con2}, we obtain that  \eqref{sistema2} is equivalent to
\begin{align}
\label{sistema3}
&2\mathrm{e}^{\frac{1}{2} \varphi}F_{\overline{\zeta}}=\dot{\varphi}\, \vert F\vert^{2}\vert G\vert^{2}\overline{G},\nonumber \\
&4\mathrm{e}^{\frac{1}{2} \varphi}G_{\overline{\zeta}}=\dot{\varphi}\, \overline{F}(1-\vert G\vert^{4})\\
& \mathrm{e}^{\frac{1}{2} \varphi}\langle\psi,\vec{e}_3\rangle_\zeta=FG.\nonumber
\end{align}
\begin{remark}
{\rm Observe that from  \eqref{sistema3}, $G$ is holomorphic if and only if $|G|\equiv 1$ and, in this case, it is clear that $G$ must be constant and  $\psi(\Sigma)$   lies on a vertical plane in $\mathbb{R}^3$.}
\end{remark}
Let us  consider $\psi_k$ a $[\varphi_k,\vec{e}_3]$- minimal surface where 
\begin{align}
& \varphi_k(z) =z, \quad z\in \mathbb{R}, \quad \text{ if $k=1$},\label{soliton}\\
&  \varphi_k(z) =\frac{2}{k-1} \log (z)  \quad z>0, \quad \text{ if $k\neq1$ }\label{singular}
\end{align}
Then, from \eqref{con1}, \eqref{con2}, \eqref{sistema3}, \eqref{soliton} and \eqref{singular}, the  Gauss map $G$ of $\psi_k$ satisfies the following complex equation 
\begin{align}
\label{equ}
&G_{\zeta\overline{\zeta}}+2\frac{\vert G\vert^{2}}{1-\vert G\vert^{4}}\overline{G}G_{\zeta}G_{\overline{\zeta}}+2k \frac{\vert G_{\overline{\zeta}}\vert^{2}}{1-\vert G\vert^{4}}G=0.
\end{align}
and
$$ \frac{G_{\overline{\zeta}}}{{1-\vert G\vert^{4}}}, \ \  \frac{\overline{G}G_{\overline{\zeta}}}{{1-\vert G\vert^{4}}}, \ \ \frac{\overline{G}^2G_{\overline{\zeta}}}{{1-\vert G\vert^{4}}}$$
are smooth functions on ${\Sigma}$. In this case, the equation \eqref{equ} gives the integrability conditions of the system \eqref{sistema3} and we have,
\begin{theorem}[{\rm \cite[Theorem 3.2]{MM2}}]
\label{W1}
Let $G$ be a not holomorphic solution of \eqref{equ} defined on a simply connected domain ${\cal U}\subset \mathbb{C}$. Then the map $\psi_1:{\cal U} \rightarrow \mathbb{R}^3$ given by
\begin{equation}
\label{rw1}
\psi_1=4\,\Re\left(\int _{\zeta_0}^\zeta\frac{\overline{G}_{\zeta}(1-G^{2})}{1-\vert G\vert^{4}} d\zeta,\int_{\zeta_{0}}^{\zeta} i\, \frac{\overline{G}_{\zeta}(1+G^{2})}{1-\vert G\vert^{4}}d\zeta, 2\int_{\zeta_{0}}^{\zeta} \frac{\overline{G}_{\zeta}G}{1-\vert G\vert^{4}}d\zeta\right)\end{equation}
is a conformal  translating soliton in $\mathbb{R}^{3}$ with Gauss map $G$.
Conversely, any translating soliton which is not on a vertical plane can be locally represented in this way.
\end{theorem}
\begin{theorem}[{\rm \cite[Theorem 3.3]{MM2}}]
\label{Wk}
Let $G$ be a not holomorphic solution of \eqref{equ} defined on a simply connected domain ${\cal U}\subset \mathbb{C}$. Then for any $k\neq 0,1$, the map $\psi_k:{\cal U} \rightarrow \mathbb{R}^3$ given by
\begin{align}
\psi_k=& \bigg( \displaystyle 4 k\Re\int_{\zeta_{0}}^{\zeta}  \frac{\overline{G}_{\zeta}(1-G^{2})}{1-\vert G\vert^{4}}\Gamma\, d\zeta, 
\displaystyle 4 k\Re\int_{\zeta_{0}}^{\zeta}  i \frac{\overline{G}_{\zeta}(1+G^{2})}{1-\vert G\vert^{4}}\Gamma\, d\zeta , 
\displaystyle\frac{2k}{k-1} \Gamma\bigg), \label{rw2}
\end{align}
where 
\begin{align*}
\Gamma = \mathrm{e}^{\displaystyle4(k-1)\Re\int_{\zeta_{0}}^{\zeta}  \frac{\overline{G}_{\zeta}G}{1-\vert G\vert^{4}}\, d\zeta}
\end{align*}
is a conformal  $\displaystyle\frac{2}{k-1}$-singular minimal surface in $\mathbb{R}^{3}_+$ with  Gauss map $G$.
Conversely, any singular minimal surface in $\mathbb{R}^3_+$  which is not on a vertical plane can be locally represented in this way.
\end{theorem}
\begin{remark} 
{\rm  We would like to point out that it came to our knowledge that Theorem \ref{W1} was also proved in \cite[Theorem 4]{Le2}. 
\\
The case $k = 0$ (i.e. of minimal surfaces in Hyperbolic space) has been studied in  \cite{K}.}
\end{remark}
\section{The Cauchy's problem.}
The Weierstrass representation described in the above Section can be applied to solve the following  general Cauchy  problem:
\begin{quote}
Let $\beta = (\beta_1,\beta_2,\beta_3):I\rightarrow \mathbb{R}^3$ be a regular analytic  curve and let $V:I\rightarrow \mathbb{S}^2$ be an analitic vector field along $\beta$ such that $\langle\beta',V\rangle=0$, $| \Pi\circ V|<1$ and $\beta_3>0$ if $k\neq1$, where $\Pi$ denotes the stereographic projection from the south pole. Find  $[\varphi_k,\vec{e}_3]$-minimal surfaces containing $\beta$ with unit normal in $\mathbb{R}^3$ along $\beta$ given by $V$.
\end{quote} 
This problem has been inspired by the classical Bj\"{o}rling problem for minimal surfaces in $\R^3$, proposed by E. G. Bj\"{o}rling in 1844  and solved by H.A. Schwarz in 1890. 
Any pair $\beta$, $V$ in the conditions of that problem {\sl a pair of Bj\"orling data}.
\begin{theorem}[{\rm \cite[Theorem 4.3]{MM2}}]
\label{cp}
For any $k\in \mathbb{R}$, $k\neq 0$, there exists a unique  $[\varphi_k,\vec{e}_3]$-minimal  which is a  solution to the Cauchy problem with  Bj\"orling data $\beta=(\beta_1,\beta_2,\beta_3)$, $V=(V_1,V_2,V_3)$. This solution, $$\psi:{\cal U}=I\times]-\epsilon,\epsilon[\subseteq\mathbb{C}\longrightarrow \mathbb{R}^3,$$ can be constructed in a neigbourhood of $\beta$ as follows: let $G:{\cal U}\rightarrow \mathbb{C}$ be the unique solution to the following system of Cauchy-Kowalewski's type, {\rm \cite{Pe}},
\begin{equation}\label{ckp}
\left\{
\begin{array}{l}\displaystyle
G_{\zeta\overline{\zeta}}+2\,\frac{\vert G\vert^{2}G_{\zeta}G_{\overline{\zeta}}}{1-\vert G\vert^{4}}\ \overline{G}+2k\,\frac{ \vert G_{\overline{\zeta}}\vert^{2}}{1-\vert G\vert^{4}}\ G=0,\\ 
\\
\displaystyle
G(s,0) =\frac{\phi_3(s)}{\phi_1(s) - i \phi_2(s)}= -\frac{\phi_1(s) + i \phi_2(s)}{\phi_3(s)}, \\ \\
G_{\overline{\zeta}}(s,0)=\left\{\displaystyle
\frac{1-|G(s,0)|^4}{4}\left(\overline{\phi}_{1}(s)+i\overline{\phi}_{2}(s) \right), \text{ if } k=1, \atop \displaystyle
\frac{1-|G(s,0)|^4}{2(k-1)\beta_{3}}\left(\overline{\phi}_{1}(s)+i\overline{\phi}_{2}(s)\right), \text{ if } k\neq 1,
\right.
\end{array}\right.
\end{equation}
where 
\begin{equation}\label{phi}
\phi(s)=(\phi_1(s),\phi_2(s),\phi_3(s))=\frac{1}{2}(\beta'(s)-i\beta'(s)\wedge V(s)), \quad u\in I.
\end{equation}
Then $\psi$ is given, up to an appropriate translation, by \eqref{rw1} if $k=1$ and  by  \eqref{rw2} if $k\neq1$ and $G$
is its  Gauss map.

\end{theorem}
\section{A Calabi's type correspondence}
Let us consider $\L^3$  the Minkowski space $\R^3$ with the Lorentz metric 
\begin{equation} \llangle \cdot ,  \cdot \rrangle = dx^2 + dy^2 - dz^2,\label{mmetric}\end{equation}
A surface in $\mathbb{L}^{3}$ is called \textit{spacelike} if the induced metric on the surface is a positive definite Riemannian metric. This kind of surfaces have played a major role in Lorentzian geometry,  for a survey of some results we refer to \cite{Bar}.

\

Let $\varphi: I\rightarrow \mathbb{R}$ be a smooth function on a real open interval $I\subseteq \mathbb{R}$. 
A spacelike surface  $ {\widetilde \Sigma}$ in $\mathfrak{D}^3=\R^2\times I\subseteq\L^3$ is called   $[\varphi, \vec{e}_3]$-{\sl maximal} if its mean curvature vector  $ {\widetilde{\text {\bf  H}}}$  satisfies
\begin{align} 
& {\widetilde {\text {\bf H}}} = \left( \overline{\nabla}^{\L^3} \varphi \right)^\perp = -\dot{\varphi} \, \vec{e}_3,\label{fmaximal}
\end{align}
where $\overline{\nabla}^{\L^3}$ denotes   the the gradient operator in $\L^3$.

As in the Euclidean case, 
a   $[\varphi, \vec{e}_3]$-maximal spacelike  surface   can be also viewed either   a critical point of the weighted volume functional 
\begin{equation}
{\widetilde V}_{\varphi} ({\widetilde \Sigma}):= \int_{{\widetilde \Sigma}}  \mathrm{e}^\varphi\ dA_{{\widetilde \Sigma}},
\end{equation}
 or a  maximal (zero mean curvature) spacelike surface in  the conformally changed metric 
\begin{equation} 
{\widetilde g}_\varphi:=  \mathrm{e}^\varphi \llangle\cdot,\cdot\rrangle.
\end{equation}

\

Well known  examples of spacelike $[\varphi, \vec{e}_3]$-maximal  are  the spacelike maximal  surfaces and the  spacelike translating solitons, whose study is an exciting and already classical mathematical research field,  see \cite{CY,QD} for some results. As in the Euclidean case, a spacelike  $[\varphi, \vec{e}_3]$-maximal   with $\varphi(p)= \alpha \log{|\llangle p,\vec{e}_3\rrangle| }$, $\alpha\in \R$, $\alpha\neq 0$ will be called {\sl singular $\alpha$-maximal} surface.

If ${\widetilde\Omega}$ is a simply connected planar domain,  the vertical graph  in $\L^3$ of a function $\tilde{u}:{\widetilde\Omega} \rightarrow \R$ is a  $[{\varphi}, \vec{e}_3]$-maximal spacelike  if and only if $\tilde{u}$ is a solution of the  following elliptic partial differential equation,
\begin{equation}
 (1-\tu_{\tx}^2) \tu_{\ty\ty} +  (1-\tu_{\ty}^2) \tu_{\tx\tx} + 2  \tu_{\tx}  \tu_{\ty} \tu_{\tx\ty} +\dot{\varphi} (\tu)\widetilde{W}^2=0, \quad (\tx,\ty)\in\wO,\label{Lfe}
\end{equation}
where  $\widetilde{W}= \sqrt{ 1-\tu_{\tx}^2- \tu_{\ty}^2}$.
The equation  \eqref{Lfe}, is equivalent to the integrability of the following differential system,
\begin{equation}  \wph_{\tx\tx}=  \frac{ 1 -   \tu_{\tx}^2}{\wW} \mathrm{e}^{\varphi(\tu)}, 
\quad 
 \wph_{\tx\ty}=  - \frac{ \tu_{\tx} \tu_{\ty}}{ \wW} \mathrm{e}^{\varphi(\tu)}, 
 \quad   
\wph_{\ty\ty}=  \frac{ 1 -   \tu_{\ty}^2}{\wW} \mathrm{e}^{\varphi(\tu)}\label{lintegra}
\end{equation}
for a convex function  $\wph:\wO \longrightarrow \R$ (unique, modulo linear polynomials).

\

In  \cite{Ca}, E. Calabi observed that there is a natural (local) connection  between Euclidean minimal graphs and Lorentzian spacelike maximal graphs which is useful  for  describing  examples and for applying similar methods to the study of their geometrical and topological properties. Recent advances of this correspondence in other ambient spaces could be found in \cite{Le,Le1,Le3,LM,PRR}. 

\

Now, we show how Calabi's correspondence  can be extended to one between the family  of $[\varphi,\vec{e}_{3}]$-minimal surfaces in $\mathbb{R}^{3}$ and the family of spacelike $[\varphi,\vec{e}_{3}]$-maximal surfaces in $\mathbb{L}^{3}$.

\begin{theorem}[{\rm \cite[Theorem 4.1]{MM1}}] Let $\Omega$ be a simply connected planar domain,  $\psi: \Omega \rightarrow \R^3$, $\psi(x,y)= (x,y,u)$ be a vertical $[\varphi ,\vec{e}_3]$-minimal  graph in $\R^3 $, $\phi$ be a solution to the system \eqref{intsystem} and $\vartheta$ be a primitive function of $\mathrm{e}^\varphi $ (that is, $\dot{\vartheta}=\mathrm{e}^\varphi)$. Then  $\widetilde{\psi}: \Omega_1 \rightarrow \L^3$ given by
\begin{equation}
 \widetilde{\psi} := (\phi_x,\phi_y, \vartheta(u) ), \label{correspondencia1}
\end{equation}
is a $[-\varphi\circ \vartheta^{-1} ,\vec{e}_3]$-maximal  spacelike graph in the Lorentz-Minkowski space whose  Gauss map
$ \widetilde{N} $ writes
\begin{align}
\widetilde{N} &= (u_x, u_y, W), 
\label{phiL3}
\end{align}
\noindent The induced metrics $g$ and $\widetilde{g}$ of $\psi$ and $\widetilde{\psi}$, respectively,  are conformal and  the mean curvature $H$ ($\widetilde{H}$) and  Gauss curvature $K$ ($\widetilde{K}$ ) of  $\psi$ ($\widetilde{\psi}$) satisfy
\begin{align}
&  \widetilde{H} + W^2 \mathrm{e}^{-\varphi(u)}   H = 0,\label{hh}\\
&\widetilde{K}+ W^4  \mathrm{e}^{-2\varphi(u)}K = 0.\label{kk}
\end{align}
\end{theorem}

\begin{remark}
Observe, that on a simply connected $[\varphi ,\vec{e}_3]$-minimal  immersion  $\psi$ in $\R^3$ the correspondence \eqref{correspondencia1} writes  as follows:
\begin{equation} \label{globalc}
\widetilde{\psi} = \int  \mathrm{e}^{\varphi (\psi_3) } (\vec{e}_3\wedge (d\psi \wedge N) +  \langle\psi,\vec{e}_3\rangle \vec{e}_3),
\end{equation}
where $\wedge$ denotes the cross product in $\R^3$. Moreover,  the singularities of $\widetilde{\psi}$ hold  where the angle  function $\langle\vec{e}_3,N\rangle$ vanishes.  \end{remark}
\begin{theorem}[{\rm \cite[Theorem 4.4]{MM1}}]
Let  $\wO$ be a simply connected planar domain, $\widetilde{\psi}: \wO \rightarrow \L^3$, $\widetilde{\psi}(\tx,\ty)= (\tx,\ty,\tu)$ be a vertical $[\varphi ,\vec{e}_3]$-maximal  graph in $\L^3 $ , $\wph$ be a solution to the system \eqref{lintegra} and $\tilde{\vartheta}$ be a primitive of $\mathrm{e}^{\varphi }$. Then the immersion given by
\begin{equation}
\psi:= (\wph_{\tx},\wph_{\ty}, \tilde{\vartheta}(\tu )),\label{correspondencia2}
\end{equation}
 is  a $[-\varphi\circ \vartheta^{-1}, \vec{e}_3]$-minimal graph  in $\R^3$   whose  induce metric, mean curvature $H$ and Gauss curvature $K$ satisfy
\begin{align}
& g:= \frac{\mathrm{e}^{ 2\varphi(\tu)}}{\wW^2} \widetilde{g},\\
&   H + \mathrm{e}^{-\varphi(\tu)} \ \wW^2\  \widetilde{H} = 0,\label{Lhh}\\
&K  + \mathrm{e}^{-2\varphi(\tu)} \ \wW^4\  \widetilde{K} = 0\label{Lkk},
\end{align}
where  $\widetilde{g}$,  $\widetilde{H}$ and $\widetilde{K}$ are the induced metric, the mean curvature and the Gauss curvature of the spacelike graph of $\tu$.
\end{theorem}
\begin{remark}
The correspondence \eqref{correspondencia2} also writes as follows
\begin{equation} \label{lglobalc}
\psi = \int  \mathrm{e}^{\varphi(\wp_3)} \left(\vec{e}_3\wedge_{\L^3} (d\widetilde{\psi}\wedge_{\L^3} \widetilde{N}) + \llangle d\widetilde{\psi},\vec{e}_3\rrangle\vec{e}_3\right),
\end{equation}
where $\wedge_{\L^3}$ denotes the cross product in $\L^3$ and $\widetilde{N}$ is the Gauss map of $\widetilde{\psi}$. The singular points of $\psi$ hold  where the angle  function $\llangle \vec{e}_3,\widetilde{N}\rrangle$ vanishes. \end{remark}
\subsubsection*{{\sc $\bullet$ Some Consequences and Applications}}
Using the above correspondence, we can prove the following statements, see \cite[Section 5]{MM1}.
\begin{enumerate}[\hspace*{0.15cm}] %
\item[-]  If  $\psi$ is a  translating soliton in $\R^3$ then $\widetilde{\psi}$ is a singular $(-1)$-maximal spacelike surface in $\L^3$, see Figure \ref{cp1}.
\item[-] If  $\widetilde{\psi}$   is a  translating soliton in $\L^3$ then $\psi$ is a singular $(-1)$-minimal  surface in $\R^3$, see Figure \ref{cp2}.
\begin{figure}[htb]
\begin{center}
\includegraphics[width=0.3\linewidth]{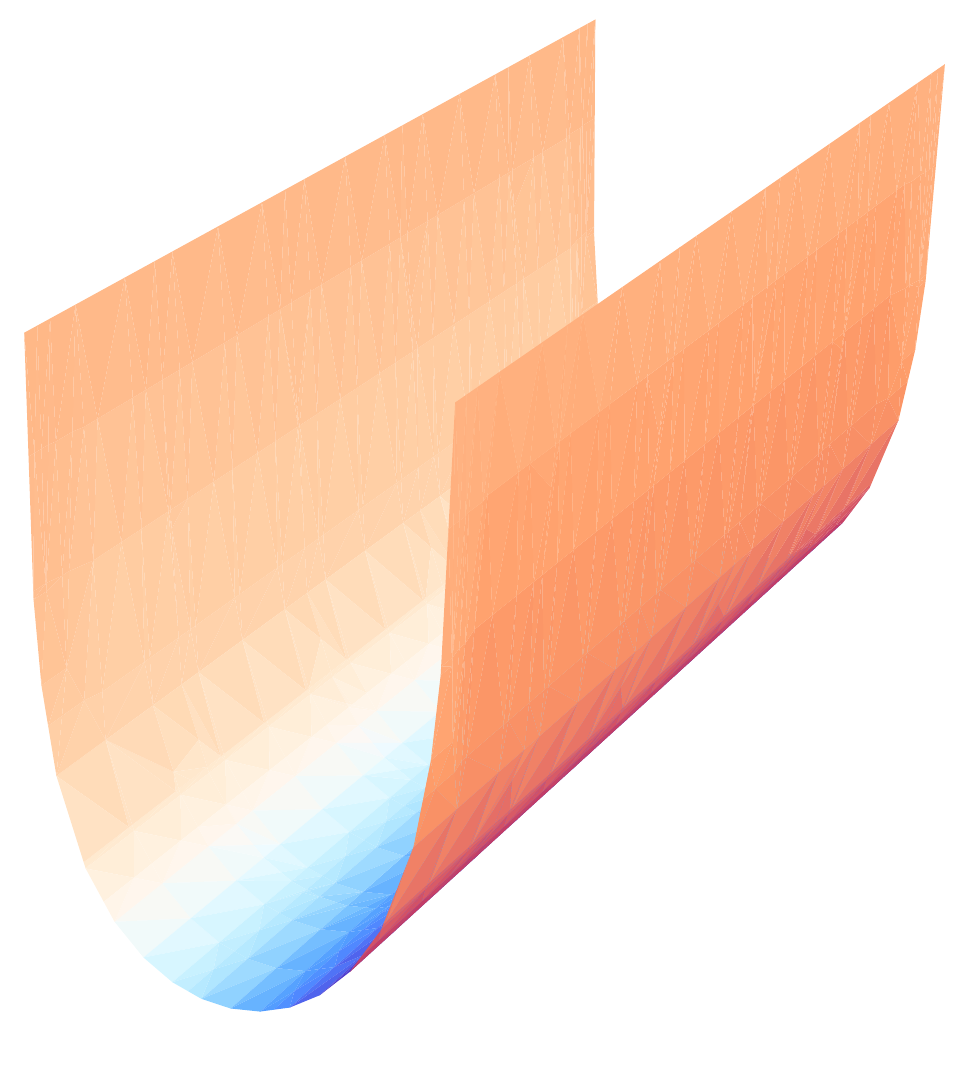}
\includegraphics[width=0.55\linewidth]{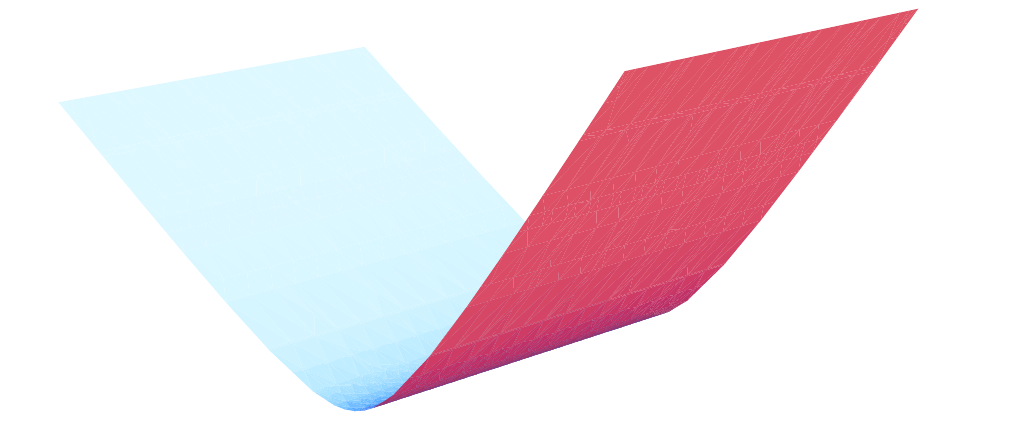}
\end{center}
\caption{Soliton  in $\R^3$ and its corresponding spacelike singular (-1)-maximal surface in $\L^3$.}\label{cp1}
\end{figure}
\begin{figure}[htb]
\begin{center}
\includegraphics[width=0.45\linewidth]{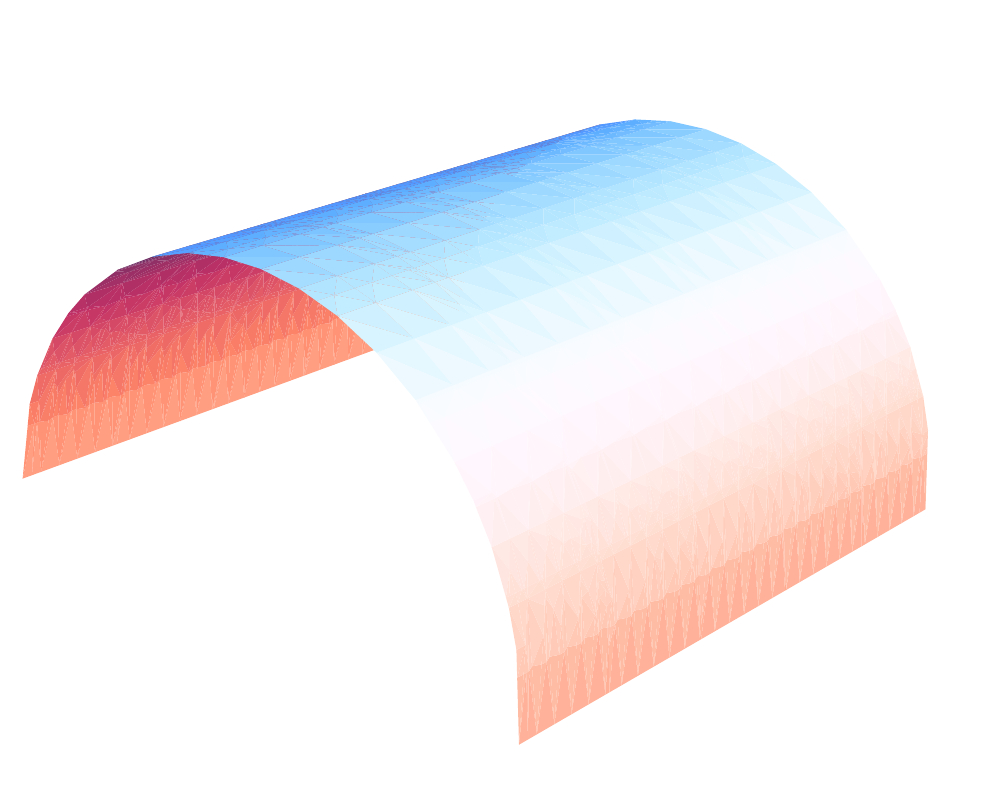}
\includegraphics[width=0.45\linewidth]{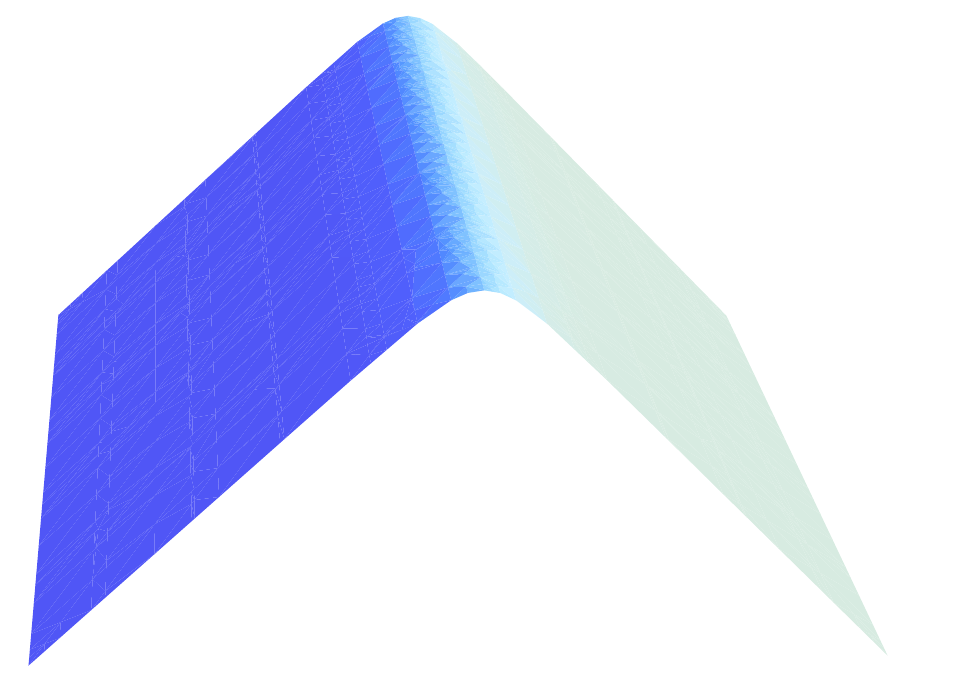}
\end{center}
\caption{singular $1$-minimal surface in $\R^3$ and  its corresponding translating soliton in $\L^3$.}\label{cp2}
\end{figure}
\item[-] Be a revolution surface with vertical rotational axis is a preserving property.
\item[-] There exists a rotationally symmetric,  entire, smooth, strictly convex spacelike translating soliton (unique up to translation) and of linear growth, Figure \ref{f4} left. 
\item[-] For any  $\alpha<-1$,  there exists a  rotationally symmetric,  entire, smooth, strictly convex spacelike singular  $\alpha$-maximal graph (unique up to  homothety) and with a linear growth, Figure \ref{f4} right.
\begin{figure}[htb] 
\begin{center}
\includegraphics[width=0.4\linewidth]{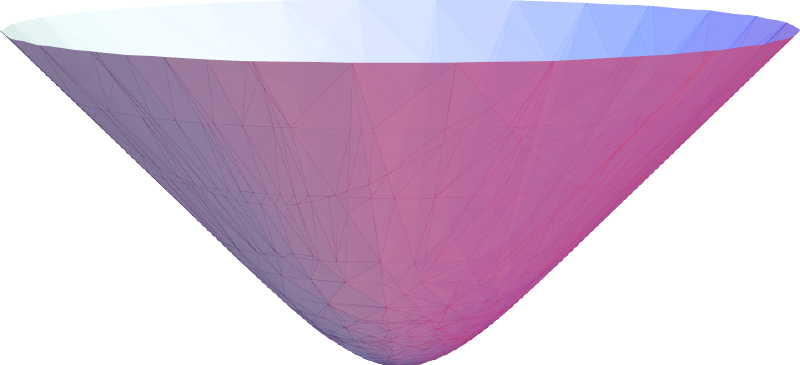}\ \ \ 
\includegraphics[width=0.4\linewidth]{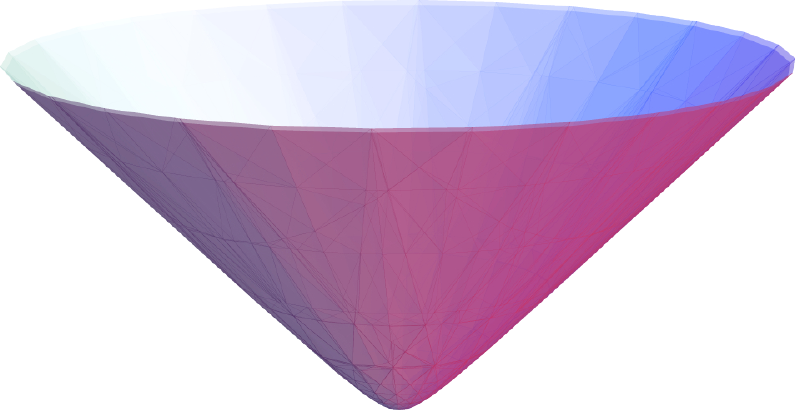}
\end{center}
\caption{Translating soliton and singular $(-2)$-maximal bowl of elliptic type in $\L^3$} \label{f4}
\end{figure}
\item[-] For any $\alpha <-1$, there exist, up to an homothety (translation), two entire spacelike singular $\alpha$-maximal  graphs (spacelike translating solitons) in $\L^3$ with linear growth and with an  isolated singularity at the origin which are asymptotic to the light cone.

\noindent {\rm This kind of examples are called either {\sl winglike solitons}  or {\sl  winglike singular $\alpha$-maximal surfaces}, respectively, see Figure \ref{fig6}}.
\begin{figure}[htb]
\begin{center}
\includegraphics[width=0.4\linewidth]{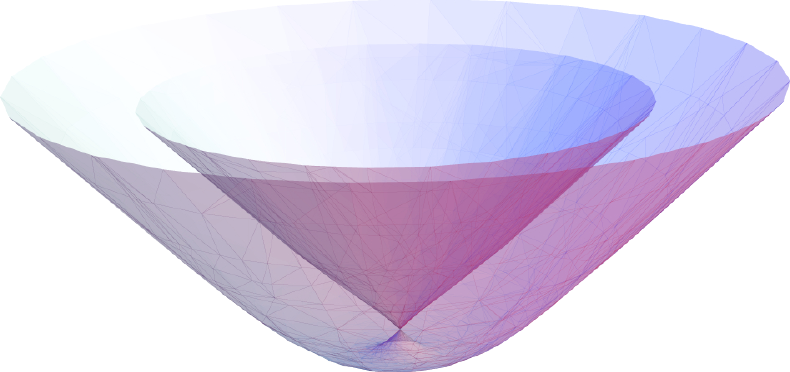} \ \ \ 
\includegraphics[width=0.4\linewidth]{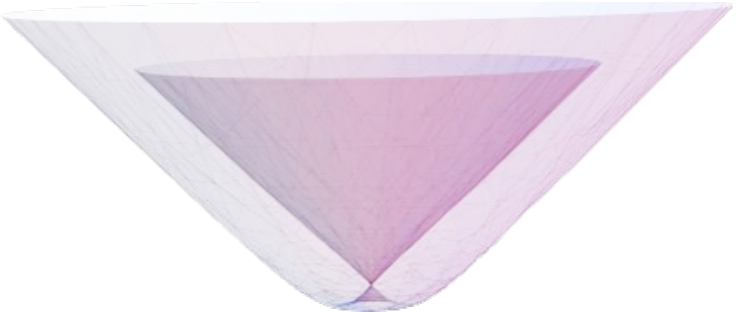}
\end{center}
\caption{Winglike soliton and singular $(-3)$-maximal winglike in $\L^3$} \label{fig6}
\end{figure}
\end{enumerate}

\end{document}